\documentclass[11pt]{article}
\usepackage{amssymb}
\usepackage{amsmath}
\usepackage{latexsym}
\newtheorem{theorem}{Theorem}[section]
\newtheorem{lemma}[theorem]{Lemma}
\newtheorem{proposition}[theorem]{Proposition}
\newtheorem{corollary}[theorem]{Corollary}

\newtheorem{remark}[theorem]{Remark}

\newcommand{\lat}{{\rm\Gamma}}
\newcommand{\dlat}{{\rm\Gamma^{\dagger}}}
\newcommand{\fd}{{\cal O}}
\newcommand{\dfd}{{{\cal O}^{\dagger}}}
\newcommand{\nl}{\newline}
\newcommand{\dist}{{\rm dist}}
\newcommand{\N}{{\bf N}}

\newcommand{\R}{{\bf R}}

\newcommand{\cF}{{\cal F}}
\newcommand{\cB}{{\cal B}}

\newcommand{\cD}{{\cal D}}

\newcommand{\cA}{{\cal A}}
\newcommand{\cV}{{\cal V}}
\newcommand{\cP}{{\cal P}}

\newcommand{\CX}{{\cal X}}
\newcommand{\bnu}{\xi}
\newcommand{\GV}{V}
\newcommand{\volume}{{\rm vol}}
\newcommand{\CB}{{\cal B}}
\newcommand{\boldeta}{\eta}
\newtheorem{thm}{Theorem}[section]
\newtheorem{cor}[thm]{Corollary}
\newtheorem{cla}[thm]{Claim}
\newtheorem{lem}[thm]{Lemma}

\newcommand{\bee}{\begin{equation}}
\newcommand{\ene}{\end{equation}}
\newcommand{\bees}{\begin{equation}}
\newcommand{\enes}{\end{equation}}
\newcommand{\bes}{\begin{split}}
\newcommand{\ens}{\end{split}}

\newcommand{\bet}{\begin{thm}}
\newcommand{\ent}{\end{thm}}
\newcommand{\bel}{\begin{lem}}
\newcommand{\enl}{\end{lem}}
\newcommand{\bec}{\begin{cor}}
\newcommand{\enc}{\end{cor}}
\newcommand{\becl}{\begin{cla}}
\newcommand{\encl}{\end{cla}}
\newcommand{\bep}{\begin{proof}}
\newcommand{\enp}{\end{proof}}
\newcommand{\ber}{\begin{rem}}
\newcommand{\enr}{\end{rem}}

\newcommand{\de}{\delta}
\newcommand{\al}{\alpha}
 \newcommand{\diam}{{\rm diam}}
 \newcommand{\dom}{{\rm Dom}}

 \newcommand{\ran}{{\rm Ran}}
 
 \newcommand{\inprod}[2]{{\langle{#1},{#2}\rangle}}
 \newcommand{\vol}{{\rm vol}}

 \newcommand{\darr}[4]{{\left\{\begin{array}{ll}
   {#1}&{#2}\\[0.2cm]
   {#3}&{#4}
 \end{array}\right.}}

\newcommand{\ia}{({\rm i})}
\newcommand{\ib}{({\rm ii})}
\newcommand{\ic}{({\rm iii})}
\newcommand{\id}{({\rm iv})}

\newcommand{\bxi}{{\bf\xi}}
\newcommand{\btheta}{{\bf\theta}}
\newcommand{\parder}[2]{\frac{\partial {#1}}{\partial {#2}}}

\title{Bethe-Sommerfeld conjecture for pseudodifferential perturbation}

\author{G. Barbatis \and L. Parnovski}

\date{ }

\begin{document}

\maketitle

\begin{abstract}
We consider a periodic pseudodifferential operator $H=(-\Delta)^l+A$ ($l>0$) in $\R^d$ which satisfies the 
following conditions: (i) the symbol of $H$ is smooth in $x$, and (ii) the perturbation 
$A$ has order smaller than $2l-1$. Under these assumptions, we prove that the spectrum
of $H$ contains a half-line. 

\

\end{abstract}

\section{Introduction}

Let $H$ be a periodic pseudodifferential elliptic operator in $\R^d$. Then its spectrum consists of
spectral bands converging to $+\infty$. These bands are separated by spectral gaps, and
one of the most important questions of the spectral theory of periodic operators is
whether the number of these gaps is finite. There is a wide belief that for $d\ge 2$, 
under fairly general conditions on $H$ the number of gaps is finite (often this
statement is called the `Bethe-Sommerfeld conjecture'). This statement (in the setting of periodic operators) is equivalent to stating that the whole interval $[\lambda,+\infty)$ is covered by the spectrum of $H$, provided $\lambda$ is big enough. If we discard certain very special cases of operators $H$ (like Schr\"odinger
operator with potential which allows the separation of variables), then, 
until recently, the conjecture was known to hold only under serious restrictions on the
dimension $d$ of the Euclidean space and the order of the operator $H$, see 
\cite{PopSkr}, \cite{Skr0}--\cite{Skr2}, \cite{Kar}, \cite{HelMoh}, \cite{PS1}, \cite{PS2}, \cite{Vel0},
\cite{Moh}.
Another type of sufficient conditions is assuming that the lattice of periods of $H$ is rational, see \cite{Skr1}, \cite{SkrSob1}, \cite{SkrSob2}.
In the paper \cite{P}, the conjecture was proved for Schr\"odinger operators with smooth periodic potentials, without any assumptions on 
dimension $d\ge 2$, or on the lattice of periods (see also \cite{Vel} for an alternative approach to this problem).  In our paper, we prove that this conjecture holds for a wide class of pseudodifferential operators.

Let us describe the results of our paper in detail. Let $h=h(x,\xi)$ be the symbol of $H$ and let $2l$ be the order of $H$ ($l>0$).
We assume a decomposition $h=h_{p}+a$, where $h_{p}\asymp|\xi|^{2l}$ (as $|\xi|\to\infty$) is the principal symbol of $H$, and $a=O(|\xi|^{\al})$ (as $|\xi|\to\infty$, where $\al<2l$) is the perturbation. The symbol $h(x,\xi)$ is assumed to be periodic in $x$ with periodicity lattice $\lat$. We denote by $\dlat$ the dual lattice, and by $\fd$ and $\dfd$ the respective fundamental domains. We also set $d(\lat)=\vol(\fd)$ and $d(\dlat)=\vol(\dfd)$.

We make the following additional assumptions:

\

(a) $h_{p}(\xi)=|\xi|^{2l}$. This assumption is made mainly for simplicity of exposition; our results are likely to hold if we replace $h_p$ by a more general principal symbol that is homogeneous in $\xi$. However it is essential that $h_p(\xi)$ does not depend on $x$ (in other words, the principal part of $H$ has constant coefficients). This latter assumption, to the best of our knowledge, is present in all approaches to proving the Bethe-Sommerfeld conjecture. 

\

(b) We assume that the symbol $a(x,\xi)$ is smooth in $x$.
This requirement is a major disadvantage of our method when compared, e.g., with the approach of Karpeshina (see \cite{Kar} and references therein).

\

(c) $a=O(|\xi|^{\al})$ and $\nabla_{\xi} a=O(|\xi|^{\al-1})$. We emphasize that we do not assume the existence of higher (than first) derivatives of $a$ with respect to $\xi$.

\

(d) Finally, we assume that $\al<2l-1$. This assumption is the most restrictive one. In particular, it
means that the results of our paper are not applicable to the Schr\"odinger operator
with a periodic magnetic potential. Indeed, in a forthcoming publication \cite{PS3} it will be shown that under conditions (a)-(c) alone, the conjecture does not hold if we only assume $\alpha<2l$. 

The method of proof in this paper follows very closely that of \cite{P}.
However, there are numerous amendments which, while being not too difficult, are not straightforward either.
We have written in detail most of the proofs, but occasionally we will refer the reader to \cite{P} if the proof of some statement in our paper is (almost) identical to the corresponding proof in \cite{P} (otherwise the size of our paper
would become almost intolerable). 
Here is the list of all the major changes we had to make to \cite{P}  to cater for a bigger class of operators: 
a generalization of the main approximation lemma has been given to be able to deal with an unbounded perturbation;
the definition of the resonance sets $\Xi_j$ has been changed; the proof of the
asymptotic formula for eigenvalues in the non-resonance region has been changed
(we use the implicit function theorem, which is slightly easier than the method used in \cite{P}); a complication arising from the fact that the mappings $F_{\xi_1,\xi_2}$
in Section 5 do not any longer provide unitary equivalence, has been addressed; finally,
Lemma \ref{lem:volinters} now includes two cases ($a$ lying inside and outside a spherical layer of radius
$2\rho$; the latter case allows a much better estimate) and the rest of Section 6 incorporates this point. Some further changes are indicated in the text. 

We have also made certain changes in order to make exposition simpler. The major 
such change is as follows:
for simplicity, we will always assume that the symbol of the perturbation $a(x,\xi)$ is a trigonometric potential in $x$, i.e. that 

\bee\label{eq:new1}
a(x,\xi)=\sum_{\theta\in\dlat,\ |\theta|<R} \hat a(\theta,\xi) e_{\theta}(x), 
\ene
where
\bee
e_{\theta}(x)=d(\lat)^{-1/2}e^{i x\cdot\theta}
\ene
and
\bee
\hat a(\theta,\xi)=\int_{\fd}a(x,\xi)e_{-\theta}(x)dx.
\ene
Here, $R$ is a fixed number. The general case of $a(x,\xi)$ being merely smooth in $x$ can be treated 
in the same way as in \cite{P}: for each large $\rho$ we choose $R=\rho^{\tau}$ with sufficiently small
(but fixed) $\tau$ and consider the truncated symbol
\bee
a'(x,\xi)=\sum_{\theta\in\dlat,\ |\theta|<R} \hat a(\theta,\xi) e_{\theta}(x).
\ene
Then the difference between the eigenvalues of the original operator $H=H_0+A$ and the truncated operator
$H'=H_0+A'$ ($A'$ is the operator with symbol $a'$) is an arbitrarily large negative power of $\rho$ (strictly speaking, we need to replace $A$ with $A'$ after our first cut-off introduced in Corollary \ref{cor:where}).
Then the results would follow if we carefully keep tracing how all the estimates for $H'$ depend on $R$. 
This has been done in detail in \cite{P}, so in order not to overburden our paper with extra notation,
we will assume that $H=H'$, i.e. that the symbol of the perturbation has the form (\ref{eq:new1}) for a fixed $R$.

\

{\bf Setting and notation.}

We fix a lattice $\lat\subset\R^d$ and denote by $\dlat$ the dual lattice.
We denote by $\fd$ and $\dfd$ the respective fundamental domains and set $d(\lat)=\vol(\fd)$, $d(\dlat)=\vol(\dfd)$.
We denote by $\cF u(\xi)$, $\xi\in\R^d$, the Fourier transform of a function $u(x)$ and by $\hat{a}(\theta)$, $\theta\in\dlat$, the Fourier
coefficients of a function $a(x)$ which is periodic with respect to the lattice $\lat$; that is
\[
(\cF u)(x)=(2\pi)^{-d/2}\int_{\R^d}u(x)e^{-ix\cdot\xi}dx \; , \qquad
\hat{a}(\theta)=d(\lat)^{-1/2}\int_{\fd}a(x)e^{-i x\cdot\theta}dx.
\]
By $\{\xi\}$ we denote the fractional part of a point $\xi\in\R^d$ with respect to the lattice $\dlat$, that is a unique point such that
$\{\xi\}\in\dfd$, $\xi-\{\xi\}\in\dlat$.
By $f\ll g$ we shall mean that there exists $0<c<\infty$ such that $f\leq cg$.

Let $r>0$. A linear subspace $V\subset\R^d$ is called a lattice subspace of dimension $n$, $1\le n \le d$, if it is spanned
by linearly independent vectors $\theta_1,\ldots,\theta_n\in\dlat$ each of which has length smaller than $r$.
We denote by $\cV(r,n)$ the set of all lattice subspaces of dimension $n$. We will usually take
$r=6MR$ for some fixed and large $M$ and $R$, so we set for simplicity $\cV(n)=\cV(6MR,n)$.
Given a subspace $V$ we denote by $\xi_V$ and $\xi_V^{\perp}$ the orthogonal projections of
$\xi\in\R^d$ on $V$ and $V^{\perp}$ respectively.
We also define
\[
\Theta_j=B(jR)\cap\dlat \; , \qquad \Theta_j'=\Theta_j\setminus\{0\}.
\]
Given $k\in\dfd$ and a set $U\subset\R^d$ we denote by $\cP^{(k)}(U)$ the orthogonal projection
in $L^2(\fd)$ onto the subspace spanned by the set $\{e_{\xi} \, : \, \{\xi\}=k \; , \;\xi\in U\}$.
Given a bounded below self-adjoint operator $T$ with discrete spectrum, we denote by $\{\mu_j(T)\}$ its eigenvalues, written
in increasing order and repeated according to multiplicity.
By $c$ or $C$ we denote a generic constant whose value may change from one line to another. However
the constants $c_1,c_2,$ etc. are fixed throughout. All constants may depend not only on the parameters of the problem,
i.e. the order $2l$, the lattice $\lat$ and the symbol $a(x,\xi)$, but also on the number $M$.

%%%%%%%%%%%%%%%%%%%%%%%%%%%%%

We consider the self-adjoint operator
\begin{equation}
H =(-\Delta)^l+A =:H_0 +A
\label{h}
\end{equation}
on $L^2(\R^d)$, where $l>0$ (not necessarily an integer)
and $A$ is a periodic PDO of order $\alpha<2l-1$
and periodicity lattice $\lat$. What we mean by this is that $A$ has the form
\[
Au(x)=(2\pi)^{-d/2}\int_{\R^d}a(x,\xi)e^{ix\cdot\xi}(\cF u)(\xi)d\xi \, ,
\]
where the symbol $a(x,\xi)$ is assumed to have the following properties: as a function of $x$ it is
$C^{\infty}$ and periodic
with periodicity lattice $\lat$; moreover, there exists $c>0$ such that
\[
|a(x,\xi)|\leq c \langle \xi\rangle^{\alpha} \; , \quad |\nabla_{\xi}a(x,\xi)|\leq c \langle \xi\rangle^{\alpha-1},
\]
for all $x,\xi\in\R^d$; here $\langle\xi\rangle =1+|\xi|$.
It is standard \cite{RS} that under the above conditions the operator $H$ admits
the Bloch-Floquet decomposition: it is unitary equivalent to a direct integral,
\begin{equation}
H\simeq\int_{\oplus}H(k)dk \; ;
\label{bf}
\end{equation}
the direct integral is taken over $\dfd$, and for each $k\in\dfd$ the operator $H(k)$ acts on
$L^2(\fd)$ as follows: It has the same symbol as $H$ and it satisfies quasi-periodic boundary conditions
depending on $k$: its domain is given by
\[
\dom(H(k)) =\{ u|_{\fd} \,  : \, u\in H^{2l}_{loc}(\R^d) \; \; , \; \; u(x+y)=e^{i k\cdot y}u(x) \mbox{ all $x\in\R^d$, $y\in\lat$}\}.
\]
When working with the operator $H(k)$ it is convenient to use the basis
$\{e_{\xi}\}_{\{\xi\}=k}\subset\dom(H(k))$. In this respect we note that
\begin{equation}
H(k)e_{\xi}=:H_0(k)e_{\xi}+A(k)e_{\xi}=
|\xi|^{2l}e_{\xi}+d(\lat)^{-1/2}\sum_{\{\eta\}=k}\hat{a}(\eta-\xi,\xi)e_{\eta}\, ,
\label{hk}
\end{equation}
for any $\xi$ with $\{\xi\}=k$.
The domain of $H(k)$ is given, equivalently, by
\[
\dom(H(k))=\{ u=\sum_{\{\xi\}=k}u_{\xi}e_{\xi} \; : \; \sum_{\{\xi\}=k}|\xi|^{4l}|u_{\xi}|^2 <\infty \} .
\]
It follows from (\ref{bf}) that $\sigma(H)=\overline{\cup_k \sigma( H(k))}$. Each $H(k)$ has  discrete spectrum,
$\sigma(H(k))=\{\lambda_j(k)\}_{j=1}^{\infty}$, where the eigenvalues are written in increasing order and
repeated according to multiplicity. By standard perturbation theory \cite{K},
each $\lambda_j(k)$ is continuous in $k\in\dfd$
and therefore for each $j$ the union $\cup_k\{\lambda_j(k)\}$ is a closed interval $[a_j,b_j]$, known as a spectral
band. It follows that $\sigma(H)=\cup_{j=1}^{\infty}[a_j,b_j]$.
In our main theorem we prove that there is only a finite number of spectral gaps. More precisely we have:
\begin{theorem}
Let $d\geq 3$. Suppose that $\rho$ is large enough. Then $\lambda=\rho^{2l}$ belongs to
$\sigma(H)$. Moreover, there exists $Z>0$ such that the interval
$[\rho^{2l}-Z\rho^{2l-d-1},\rho^{2l}+Z\rho^{2l-d-1}]$
lies inside a single spectral band.
\label{thm:bs}
\end{theorem}
A similar statement is valid for $d=2$; see Theorem \ref{thm:rfin}.

\section{Preliminary results}

\subsection{Abstract results}

In this subsection we present some abstract theorems about
self-adjoint operators with discrete spectra. The following two
lemmas have been proved in \cite{P}.  Roughly speaking, they state that that
the eigenvalues of a perturbed operator $H=H_0+A$ that lie in a
specific interval $J$ are very close to those of $\sum_P PHP$,
where $\{P\}$ is a carefully chosen family of eigenprojections of $H_0$. 
\begin{lemma}
Let $H_0$, $A$ and $B$ be self-adjoint operators such that $H_0$
is bounded below and has compact resolvent, and $A$ and $B$ are
bounded. Put $H=H_0+A$ and $\hat H=H_0+A+B$ and denote by
$\mu_l=\mu_l(H)$ and $\hat\mu_l=\mu_l(\hat H)$ the
eigenvalues of these operators. Let $\{P_j\}$ ($j=0,\dots,n$) be a
collection of orthogonal projections commuting with $H_0$ such
that $\sum P_j=I$, $P_jAP_k=0$ for $|j-k|>1$, and $B=P_nB$. Let
$l$ be a fixed number. Denote by $a_j$ the distance from $\mu_l$
to the spectrum of $P_jH_0P_j$.
Assume that for $j\ge 1$ we have $a_j>4a$, where $a:=\|A\|+\|B\|$. Then
$|\hat\mu_l-\mu_l|\le 2^{2n}a^{2n+1}\prod_{j=1}^n(a_j-2a)^{-2}$.
\label{ll1}
\end{lemma}
{\em Proof.} This is Lemma 3.1 of \cite{P}.

\begin{lemma}\label{new:1}
Let $H_0$ and $A$ be self-adjoint
operators such that $H_0$ is bounded below and has compact
resolvent and $A$ is bounded. Let $\{P^m\}$ ($m=0,\dots,n$) be a
collection of orthogonal projections commuting with $H_0$ such
that if $m\ne n$ then $P^mP^n=P^mAP^n=0$. Denote $Q:=I-\sum P^m$.
Suppose that each $P^m$ is a further sum of orthogonal projections
commuting with $H_0$: $P^m=\sum_{j=0}^{j_m} P^m_j$ such that
$P_j^mAP_l^m=0$ for $|j-l|>1$ and $P_j^mAQ=0$ if $j<j_m$. Let
$b:=\|A\|$ and let us fix an interval $J=[\lambda_1,\lambda_2]$ on the
spectral axis which satisfies the following properties: the spectra of
the operators $QH_0Q$ and $P_j^kH_0P_j^k$, $j\ge 1$ lie outside $J$;
moreover, the distance from the spectrum of $QH_0Q$ to $J$ is
greater than $6b$ and the distance from the spectrum of
$P_j^kH_0P_j^k$ ($j\ge 1$) to $J$, which we denote by $a_j^k$, is
greater than $16b$. Denote by $\mu_p\le\dots\le\mu_q$ all
eigenvalues of $H:=H_0+A$ which are inside $J$. Then the
corresponding eigenvalues $\tilde\mu_p,\dots,\tilde\mu_q$ of the
operator
\[
\tilde H:=\sum_m P^mHP^m+QH_0Q
\]
are eigenvalues of
$\sum_m P^mHP^m$, and they satisfy
\[
|\tilde\mu_r-\mu_r|\le
\max_m \left[(6b)^{2j_m+1}\prod_{j=1}^{j_m}(a_j^m-6b)^{-2}\right];
\]
all other
eigenvalues of $\tilde H$ are outside the interval
$[\lambda_1+2b,\lambda_2-2b]$.

Moreover, there exists an injection $G$ defined on the set of eigenvalues of the operator $\sum_m P^mHP^m$ (all eigenvalues are counted according to their multiplicities) and mapping them to the
set of eigenvalues of $H$ (again considered counting multiplicities) such that:

$\ia$ all eigenvalues of $H$ inside $[\lambda_1+2b,\lambda_2-2b]$ have a pre-image,

$\ib$ if $\mu_j\in [\lambda_1+2b,\lambda_2-2b]$ is an eigenvalue of $\sum_m P^mHP^m$, then
\[
|G(\mu_j)-\mu_j|\leq
\max_m \left[(6b)^{2j_m+1}\prod_{j=1}^{j_m}(a_j^m-6b)^{-2}\right],
\]

and

$\ic$ we have $G(\mu_j(\sum_m P^mHP^m))=\mu_{j+l}(H)$, where $l$ is the number of eigenvalues of
$QH_0Q$ which are smaller than $\lambda_1$.

\label{ll2}
\end{lemma}
{\em Proof.} This is Lemma 3.2 and Corollary 3.3 of \cite{P}.

The last two lemmas involve bounded perturbations. The next
proposition 
is a generalization of Lemma \ref{new:1} to the case where the perturbation is unbounded. 
\begin{proposition}
Let $H_0$ and $A$ be self-adjoint operators such that $H_0$ and $H=H_0+A$ are bounded below and have
compact resolvents. Assume that $A$ is bounded relative to $H_0$ and that there exists
$\epsilon \in (0,1)$ and $k_{\epsilon}>0$ such that
\begin{equation}
|\inprod{Au}{u}|\leq \epsilon\inprod{H_0u}{u}+k_{\epsilon}\|u\|^2, \;\; u\in\dom(H_0).
\label{margo}
\end{equation}
Let $P_0,\ldots,P_N$ be orthogonal projections commuting with $H_0$ such that $P_iAP_j=0$ if $|i-j|>1$.
Let $P=\sum_{i=0}^NP_i$, $Q=I-P$, and assume that $P_iAQ=0$ if $i<N$. Assume that $AP$ is bounded and
let $b=\|AP\|$. Let $J=[\lambda_1,\lambda_2]$ be an interval on the spectral axis such that
\begin{equation}
\dist(\sigma(QH_0Q),J)\geq D_1, \quad \dist(\sigma(P_jH_0P_j),J)\geq D_2 \, , \;\; j\geq 1,
\label{dio}
\end{equation}
where $D_1:=(4b+2(\epsilon\lambda_2+k_{\epsilon}))/(1-\epsilon)$ and $D_2:=20b$.
Let $\tilde{H}=PHP +QHQ$. Then the following holds true: given an eigenvalue $\mu_l(H)$ of $H$ inside $J$, the corresponding
eigenvalue $\mu_l(\tilde{H})$ of $\tilde{H}$ is an eigenvalue of $PHP$. Moreover,
$|\mu_l(H)-\mu_l(\tilde{H})|<3b/4^N$.
\label{prop:abstract}
\end{proposition}
{\em Proof.} It follows from our assumptions that $H=\tilde{H}+P_NAQ +QAP_N$. Therefore,
\begin{equation}
\tilde{H} -2b(P_N+Q)\leq H \leq \tilde{H} + 2b(P_N+Q),
\label{ypou}
\end{equation}
and, in particular,
\begin{equation}
\mu_l(\tilde{H} -2b(P_N+Q))\leq \mu_l(H) \leq \mu_l(\tilde{H} + 2b(P_N+Q)).
\label{ypou1}
\end{equation}
Notice that the operator $\tilde{H}\pm 2b(P_N+Q)$ can be decomposed as $(PHP\pm 2bP_N)\oplus (QHQ \pm 2bQ)$.

{\bf Claim.} $\mu_l(\tilde{H}+ 2b(P_N+Q))$ is an eigenvalue of $PHP+ 2bP_N$. Indeed, suppose it is an eigenvalue
of $QHQ + 2bQ$, say $\mu_l(\tilde{H}+ 2b(P_N+Q))=\mu_i(QHQ)+ 2b$. Then the fact that $\mu_l(H)\in J$ implies that
\begin{equation}
\lambda_1-2b \leq\mu_i(QHQ)\leq\lambda_2+2b\, .
\label{cor}
\end{equation}
Moreover, from (\ref{margo}) and min-max we have
\begin{equation}
(1-\epsilon)\mu_i(QH_0Q)-k_{\epsilon}\leq \mu_i(QHQ)\leq (1+\epsilon)\mu_i(QH_0Q)+k_{\epsilon}.
\label{cor1}
\end{equation}
Now, we have either $\mu_i(QH_0Q)\leq\lambda_1-D_1$ or $\mu_i(QH_0Q)\geq\lambda_2+D_1$. In the first case (\ref{cor}) and (\ref{cor1})
give $\lambda_1-2b\leq (1+\epsilon)(\lambda_1-D_1)+k_{\epsilon}$, which contradicts the definition of $D_1$.
In the second case we obtain $(1-\epsilon)(\lambda_2+D_1)-k_{\epsilon}\leq \lambda_2+2b$, which is also a contradiction.
Hence the claim has been proved.

So $\mu_l(\tilde{H}+ 2b(P_N+Q))$ is an eigenvalue of $PHP+ 2bP_N$,
$\mu_l(\tilde{H}+ 2b(P_N+Q))=\mu_i(PHP+2bP_N)$, say. Let
\[
a_j=\dist(\mu_i(PHP),\sigma(P_jH_0P_j)) \; , \quad 1\leq j\leq N.
\]
From (\ref{ypou}) we have $\mu_i(PHP)\in
[\lambda_1-2b,\lambda_2+2b]$. Hence for $j\geq 1$,
$a_j\geq D_2-2b=18b$.
We can now apply Lemma \ref{ll1} to the unperturbed operator $PHP=PH_0P+PAP$ with the perturbation $B=2bP_N$. 
We conclude that
\[
|\mu_i(PHP)-\mu_i(PHP+2bP_N)| \leq 2^{2N}(3b)^{2N+1}\prod_{j=1}^N(a_j-6b)^{-2} \leq \frac{3b}{4^N},
\]
completing the proof of the proposition. $\hfill \Box$

%%%%%%%%%%%%%%
\subsection{Perturbation cut-off}
%%%%%%%%%%%%%%

We now return to the operator $H$ introduced in (\ref{h}).
We fix $\rho>0$; our aim is to prove that if $\rho$ is large enough then $\rho^{2l}\in\sigma(H)$.
Many of the statements that follow are valid provided $\rho$ is sufficiently large; usually this will not be mentioned explicitly.

\begin{lemma}\label{lem:new2}
Let $k\in\dfd$ and $U,V\subset\R^d$ be such that $\dist(U,V)>R$. Then
\begin{equation}
\cP^{(k)}(V)A(k)\cP^{(k)}(U)=0
\label{pvp}
\end{equation}
\label{lem:pvp}
\end{lemma}
{\em Proof.} We simply note that, because of (\ref{hk}), if
$\xi\in U$, $\{\xi\}=k$, then $A(k)e_{\xi}$ is a linear
combination of $e_{\eta}$, $\{\eta\}=k$, $\eta \in \xi+ B(R)$.
$\hfill \Box$

Let $\chi=\chi_{\rho}$ be the characteristic function of the
set $\{ \xi\in\R^d \; : \;  ||\xi| -\rho|<\rho/2\}$. Define the projection
$\bar{P}=\cF^{-1}\chi\cF$. We have
\begin{lemma}
There exists $L>0$ such that $\|A\bar{P}\|\leq L \rho^{\alpha}$. 
\label{lem:norm}
\end{lemma}
{\em Proof.}
Let $a_{\bar{P}}(x,\xi):=a(x,\xi)\chi(\xi)$ be the symbol of $A\bar{P}$. Then $A\bar{P}=\sum_{\theta\in\dlat}A_{\theta}$, where
\[
A_{\theta}u(x)= (2\pi)^{-d/2}d(\lat)^{-1/2}\int_{\R^d}\hat{a}_{\bar{P}}(\theta,\xi)e^{i(\theta+\xi)\cdot x}(\cF u)(\xi)d\xi \,.
\]
The smoothness of $a(x,\xi)$ with respect to $x$ implies that $|
\hat{a}(\theta,\xi)| \ll\langle\theta\rangle^{-N}\langle\xi\rangle^{\alpha}$ for any
$N\in\N$; hence
\[
| \hat{a}_{\bar{P}}(\theta,\xi)| \ll
\langle\theta\rangle^{-N}\rho^{\alpha} \; , \;\; N\in\N \, .
\]
It follows that for any $u\in L^2(\R^d)$,
\begin{eqnarray*}
\|A_{\theta}u\|^2&=&(2\pi)^{-d}d(\lat)^{-1}\int_{\R^d} \bigg| \int_{\R^d}\hat{a}_{\bar{P}}(\theta,\xi)e^{i\theta\cdot x}e^{ i\xi\cdot x}
(\cF u)(\xi)d\xi\bigg|^2 dx\\
&=&d(\lat)^{-1}\int_{\R^d} | \hat{a}_{\bar{P}}(\theta,\xi)(\cF u)(\xi)|^2 d\xi\\
&\leq& c\langle\theta\rangle^{-2N}\rho^{2\alpha}\|u\|^2 \,.
\end{eqnarray*}
Taking $N>d$, we conclude that
\[
\|A\bar{P}\| \leq \sum_{\theta\in\dlat}\|A_{\theta}\| \ll\rho^{\alpha}\sum_{\theta}\langle\theta\rangle^{-N}
\ll\rho^{\alpha},
\]
as required. $\hfill \Box$

Let $\bar{P}(k)$ be the orthogonal projection on the linear span
of the set $\{e_{\xi}\; : \; \{\xi\}=k \; , \;
||\xi|-\rho|<\rho/2\}$. It is easily verified that the Floquet
decomposition of $A\bar{P}$ is $\int_{\oplus}A(k)\bar{P}(k)dk$,
that is $(A\bar{P})(k)=A(k)\bar{P}(k)$. This implies in particular
that
\begin{equation}
\|A(k)\bar{P}(k)\|\leq\|A\bar{P}\| \leq L\rho^{\alpha}.
\label{kosm}
\end{equation}
Let $L>0$ be as in Lemma \ref{lem:norm} and put $J=[\rho^{2l} -100 L\rho^{\alpha},\rho^{2l}+100 L\rho^{\alpha}]$.
\begin{corollary}
Let $k\in\dfd$ be fixed. Let $\bar{P}=\bar{P}(k)$ be the
orthogonal projection on the linear span of $\{e_{\xi}\; : \;
\{\xi\}=k \; , \;  ||\xi|-\rho|<\rho/2\}$ and let $Q=I-\bar{P}$.
Define $\tilde{H}(k):=\bar{P}H(k)\bar{P}$. Then the
following holds true for $\rho$ large enough: given an eigenvalue
$\mu_l(H(k))$ of $H(k)$ inside $J$, the corresponding eigenvalue
$\mu_l(\bar{P}H(k)\bar{P}+QH(k)Q)$ is an eigenvalue of $\tilde H(k)$;
moreover, there exists $c>0$, independent of $k\in\dfd$, such that
$| \mu_l(H(k))-\mu_l(\bar{P}H(k)\bar{P}+QH(k)Q)| < \exp(-c\rho)$.
\label{cor:where}
\end{corollary}
{\em Proof.} We shall apply Proposition \ref{prop:abstract} to the operator $H(k)=H_0(k)+A(k)$. We fix
a natural number $N$ (to be determined later) and we write $\bar{P}=\oplus_{j=0}^NP_j$, where $P_0$ is the orthogonal projection on the liner span of
\[
\{ e_{\xi} \; : \; \{\xi\}=k \; , \;  | |\xi|-\rho|<\rho/4\},
\]
and $P_j$, $j\geq 1$, is similarly defined for
\[
\{ e_{\xi} \; : \; \{\xi\}=k \; , \; \frac{\rho}{4}+\frac{(j-1)\rho}{4N}\leq | |\xi|-\rho|<\frac{\rho}{4}+\frac{j\rho}{4N}\},
\]
It follows from Lemma \ref{lem:pvp} that
$P_iA(k)P_j=0$ if $|i-j|>1$ and, similarly, $P_jA(k)Q=0$ for $j<N$. We also note that 
\bee
b:=\|A(k)\bar{P}\|\leq L\rho^{\alpha}
\ene
by (\ref{kosm}).

The inequality $s^{\alpha}<s^{2l}+1$, after substituting $s=|\xi|\epsilon^{\frac{1}{2l-\al}}$,  implies that $|\xi|^{\alpha}\leq \epsilon|\xi|^{2l}+\epsilon^{-\alpha/(2l-\alpha)}$
for any $\epsilon>0$.
Thus, for any $\xi$ with $\{\xi\}=k$,
\[
\Big|\sum_{\{\eta\}=k}\hat{a}(\xi-\eta,\eta)\Big| \leq \epsilon |\xi|^{2l} +
c\epsilon^{-\frac{\alpha}{2l-\alpha}}.
\]
This implies that (\ref{margo}) is valid with $k_{\epsilon}=c
\epsilon^{-\frac{\alpha}{2l-\alpha}}$.
Hence $D_1\ll \rho^{\alpha}$. Since \nl $\dist(\sigma(QH_0Q),J)\geq c\rho^{2l}$, the first relation in (\ref{dio}) is satisfied;
so is the second by a similar argument.

So all assumptions of Proposition \ref{prop:abstract} are fulfilled. We conclude that $|\mu_l(H(k))-\mu_l(\tilde{H}(k))|<3b/4^N$.
Taking $N =[\rho]-1$ concludes the proof. 
$\hfill \Box$

This corollary shows that we can study the spectrum of $\tilde H(k)$ instead
of $H(k)$.

%%%%%%%%%%%%%%%%%%%%%%%%%%%%%%%%%%%%%%%%%

\section{Reduction to invariant subspaces}

%%%%%%%%%%%%%%%%%%%%%%%%%%%%%%%%%%%%%%%%%

The Floquet decomposition and Corollary \ref{cor:where} has led us
to the study of the eigenvalues of $\tilde H(k)$,
$k\in\dfd$, that are close to $\lambda=\rho^{2l}$. 
The operator $\tilde H(k)$ which was defined in Corollary \ref{cor:where} is a bounded perturbation of $\bar{P}H_0(k)\bar{P}$ (for a fixed $\rho$), so we can
apply Lemma \ref{ll2} to it. This will
require a specific choice of the projections $\{P^k_j\}$.
Because they
have to be invariant for $H_0(k)$, they will be of the form
$\cP^{(k)}(U)$ for some carefully defined sets $U\subset\R^d$,
localized near $|\xi|=\rho$.

We define the spherical layer
\[
\cA=\{\bxi\in\R^d \; : \; | |\bxi|^{2l}-\rho^{2l}|<100 L\rho^{\alpha} \}.
\]
It has width of order $\rho^{\alpha-2l+1}$.
Note that for all $\xi\in\cA$ we have
\begin{equation}
| |\xi|^2-\rho^2 |\ll\rho^{\alpha-2l+2}.
\label{zeta1}
\end{equation}
We fix numbers $q_0,q_1,\ldots,q_{d-1}$ and $\gamma$ such that
\begin{equation}
0<q_0<q_1<\ldots<q_{d-1}<1 \qquad , \qquad
\alpha<\gamma < 2l-2+q_0 \;  ,
\label{assumptions}
\end{equation}
and also
\begin{equation}
q_1\geq \frac{3+\alpha-2l}{2}.
\label{q1}
\end{equation}
The existence of such numbers follows from our assumption $\alpha<2l-1$. 
We also define
\[
\epsilon_0:=\frac{1}{100}\min\{  2l-2+q_0-\alpha, \gamma-\alpha, 1-q_{d-1}, \min_s(q_s-q_{s-1})\}
\]
Given a lattice subspace $V\in\cV(n)$, we define the sets
\begin{eqnarray*}
&& \Xi_0(V)=\{ \xi\in\cA \; : \;  |\bxi_V|<\rho^{q_n}\}; \\
&& \Xi_1(V)= (\Xi_0(V) + V)\cap\cA \; ; \\
&& \Xi_2(V)=\Xi_1(V)\setminus \Big( \bigcup_{m=n+1}^{d-1}\bigcup_{W\in\cV(m), W\supset V}\Xi_1(W)\Big)\; ; \\
&& \Xi_3(V)=\{ \xi=\eta +\btheta \; : \; \eta\in \Xi_2(V) \; , \; \theta\in V\cap \dlat \; , \;
| |\xi|^{2l}-\rho^{2l}|<\rho^{\gamma}\}\, ; \\
&&\Xi(V)=\Xi_3(V)+\Theta_M\; .
\end{eqnarray*}
We also define
\[
\cD=\bigcup_{m=1}^{d-1}\bigcup_{W\in\cV(m)}\Xi_2(W), \qquad \cB=\cA\setminus\cD .
\]
Note that $\Xi_2(V)\subset\Xi_3(V)$ by (\ref{assumptions}) and $\Xi_3(\{0\})=\cB$.
We also have
\begin{equation}
| |\xi|^2-\rho^2 |<c\rho^{\gamma-2l+2}
\label{bo}
\end{equation}
for all $\xi\in\Xi_3(V)$.

We now proceed to establish further properties of these sets. Let us
stress again that in what follows we shall often implicitly assume
that $\rho$ is sufficiently large.
\begin{lemma}
Let $V\in\cV(n)$, $0\leq n\leq d-1$. Then for $\xi\in\Xi_1(V)$ we have
$|\xi_V| <2\rho^{q_n}$.
\label{lem:1}
\end{lemma}
{\em Proof.} Let $\xi'\in\Xi_0(V)$ be such that $\xi-\xi'\in V$. Then
\[
|\xi_V|^2 \leq | |\xi_V|^2-|\xi'_V|^2| + |\xi'_V|^2 =| |\xi|^2-|\xi'|^2|+|\xi'_V|^2
<c\rho^{\alpha-2l+2}+\rho^{2q_n} \; ,
\]
from which the result follows. $\hfill \Box$

\begin{lemma}
Let $V\in\cV(n)$, $0\leq n\leq d-1$. Then for $\xi\in\Xi_3(V)$ we have
$|\xi_V| <2\rho^{q_n}$.
\label{ba}
\end{lemma}
{\em Proof.} Let us write $\xi=\eta+\theta$ with $\eta\in \Xi_2(V)$, $\theta\in V\cap\dlat$. Then, by (\ref{zeta1}),
\[
|\xi_V^{\perp}|^2=|\eta|^2-|\eta_V|^2\geq \rho^2-c\rho^{\alpha-2l+2}-\rho^{2q_n},
\]
and therefore, by (\ref{bo}),
\[
|\xi_V|^2=|\xi|^2-|\xi_V^{\perp}|^2\leq (\rho^2 +c\rho^{\gamma-2l+2}) -(\rho^2-c\rho^{\alpha-2l+2}-\rho^{2q_n})
\leq 2\rho^{2q_n},
\]
and the result follows.  $\hfill \Box$

\begin{corollary}
If $\xi\in\Xi(V)$, then $|\xi_V|\ll\rho^{q_n}\,$.
\label{cor:1}
\end{corollary}

\begin{lemma}
Let $V\in\cV(n)$, $0\leq n\leq d-1$.
Let $\xi\in\Xi_3(V)$, $\xi=\eta+\theta$ with $\eta\in\Xi_2(V)$ and $\theta\in V\cap\dlat$. Then
$|\theta|\ll\rho^{q_n}$.
\label{lem:estth}
\end{lemma}
{\em Proof.} We have $|\theta|=|\theta_V|\leq|\xi_V|+|\eta_V|\leq c\rho^{q_n}$. $\hfill \Box$

\begin{lemma}
If $\xi\in\Xi_3(V)$ then
\begin{equation}
\ia \;\; |  |\xi| -\rho | <c\rho^{\gamma-2l+1} \quad , \quad \ib \;\; |   |\xi_V^{\perp}| -\rho |<c\rho^{2q_n-1}\; .
\label{paok}
\end{equation}
\end{lemma}
{\em Proof.} Part $\ia$ follows directly from the definition of $\Xi_3(V)$. Part $\ib$ follows from $\ia$ and Lemma \ref{ba}. $\hfill \Box$

The following geometric lemma will be used repeatedly in what follows.
\begin{lemma} We have
\begin{equation}
|\xi_{V_1+V_2}|\ll |\xi_{V_1}|+|\xi_{V_2}|\, ,
\label{v1v2}
\end{equation}
for any two subspaces $V_1$ and $V_2$ generated by vectors in $\dlat\cap B(R)$.
\label{lem:geom}
\end{lemma}
{\em Proof.} Suppose that $V_1$ and $V_2$ are two lattice subspaces such that $W\ne V_1$ and $W\ne V_2$, where
$W:=V_1+V_2$ (otherwise the estimate is trivial). Let $\phi$ be the angle between $V_1$ and $V_2$. 
This means that $\phi$ is a minimum of angles between $\xi_1$ and $\xi_2$ where $\xi_j\in U_j:=W\ominus V_j$
(the orthogonal complement). Then a simple geometry shows that 
\bee
|\xi_{V_1+V_2}|\le \frac{|\xi_{V_1}|+|\xi_{V_2}|}{\sin(\phi/2)}.
\ene
Since the number of pairs $(V_1,V_2)$ is finite, this proves the statement. 
$\hfill\Box$

\begin{remark}
The key part of extending the proof to the case when the symbol $a$ is just smooth 
in $x$ (and is no longer a trigonometric polynomial in $x$) is checking how
the constant in (\ref{v1v2}) depends on $R$. This has been done in detail
in Section 4 of \cite{P}.
\end{remark}

\begin{lemma}
Let $V_i\in\cV(n_i)$, $i=1,2$, be two lattice subspaces such that neither of them is contained in the other and assume that $n_2\geq n_1$.
Then for any $\xi_i\in\Xi_2(V_i)$, $i=1,2$, we have:
$|\xi_1-\xi_2|>\rho^{q_{n_2}+\epsilon_0}$.
\label{lem:wind}
\end{lemma}
{\em Proof.}
Suppose to the contrary that $|\xi_1-\xi_2|\leq\rho^{q_{n_2}+\epsilon_0}$.  From Lemma \ref{lem:1} we have
$|(\xi_i)_{V_i}|\leq 2\rho^{q_{n_i}}$ and therefore
\[
|(\xi_1)_{V_2}|\leq |\xi_1-\xi_2|+|(\xi_2)_{V_2}|\leq \rho^{q_{n_2}+\epsilon_0}+ 2\rho^{q_{n_2}}\leq 2\rho^{q_{n_2}+\epsilon_0}\; .
\]
Letting $W=V_1+V_2$ and $m=\dim(W)$ we have $m>n_2$ and hence, by (\ref{v1v2}),
\[
|(\xi_1)_W| <c(|(\xi_1)_{V_1}|+|(\xi_1)_{V_2}|) <c(\rho^{q_{n_1}}+ 2\rho^{q_{n_2}+\epsilon_0})<\rho^{q_m}.
\]
Hence $\xi_1\in\Xi_1(W)$, which is a contradiction. $\hfill \Box$

\begin{lemma}
Let $V_i\in\cV(n_i)$, $i=1,2$, be two lattice subspaces such that neither of them is contained in the other.
Then for any $\xi_i\in\Xi_3(V_i)$, $i=1,2$, there holds
$|\xi_1-\xi_2|>\max\{\rho^{q_{n_1}},\rho^{q_{n_2}}\}$.
\label{lem:wind1}
\end{lemma}
{\em Proof.} Assume that $n_2\geq n_1$. Writing
$\xi_i=\eta_i+\theta_i$ with $\eta_i\in\Xi_2(V_i)$ and $\theta_i\in V_i\cap\dlat$, we have from Lemmas \ref{lem:estth} and \ref{lem:wind},
\[
|\xi_1-\xi_2|\geq |\eta_1-\eta_2|-|\theta_1|-|\theta_2| \geq \rho^{q_{n_2}+\epsilon_0}-c\rho^{q_{n_2}} \geq \rho^{q_{n_2}}\, .
\]

\begin{proposition}
If $V_i\in\cV(n_i)$, $i=1,2$, are two different lattice subspaces,
then $(\Xi(V_1)+\Theta_M)\cap(\Xi(V_2)+\Theta_M)=\emptyset$.
\label{prop:1}
\end{proposition}
{\em Proof.} If neither of $V_1$, $V_2$ is contained in the other, then the result follows from Lemma \ref{lem:wind1},
so we assume that $V_1\subset V_2$.
We consider $\xi_i\in\Xi_3(V_i)$ and shall prove that the difference $\theta=\xi_1-\xi_2$ cannot belong
to $\Theta_{4M}$. Let $\eta_i\in\Xi_2(V_i)$ and $\theta_i\in V_i\cap\dlat$ be such that $\xi_i=\eta_i+\theta_i$.
We distinguish two cases: \nl
$\ia$ $\theta\in V_2$. In this case we write $\eta_2=\tilde{\eta}+a$ with $\tilde{\eta}\in\Xi_0(V_2)$ and $a\in V_2$.
We then obtain $\eta_1=\tilde{\eta}+(a+\theta_2+\theta-\theta_1)\in\Xi_1(V_2)$, which contradicts the fact that
$\eta_1\in\Xi_2(V_1)$. \nl
$\ib$ $\theta\not\in V_2$. We argue again by contradiction, assuming that $\theta\in\Theta_{4M}$.
Then, in particular, $|\theta|\gg 1$.
We claim that
$|\eta_2\cdot\theta|>\rho^{q_{n_2}+\epsilon_0}|\theta|$; indeed, if this were not the case then we would have,
with $U$ being the linear span of $V_2$ and $\theta$ (and thus a lattice subspace):
\begin{eqnarray*}
|(\eta_2)_{U}|&\leq& c(|(\eta_2)_{V_2}|+|\eta_2\cdot\theta|) \\
&\leq&c (\rho^{q_{n_2}}+\rho^{q_{n_2}+\epsilon_0})\\
&\leq&\rho^{q_{n_2+1}}.
\end{eqnarray*}
Therefore, $\eta_2\in \Xi_1(U)$, which is a contradiction. Hence
\[
|\xi_2\cdot\theta|\geq |\eta_2\cdot\theta| -|\theta_2\cdot\theta| \geq \rho^{q_{n_2}+\epsilon_0}
-12R\rho^{q_{n_2}}\geq\rho^{q_{n_2}+\epsilon_0}/2\; ,
\]
and therefore
\begin{eqnarray*}
| |\xi_1|^{2l}-\rho^{2l} |&\geq& | |\xi_2+\theta|^{2l}-|\xi_2|^{2l} |- | |\xi_2|^{2l}-\rho^{2l} |\\
&\geq& c\rho^{2l-2+q_{n_2}+\epsilon_0}   -\rho^{\gamma} \\
&\geq&\rho^{\gamma} \, ,
\end{eqnarray*}
which is a contradiction. $\hfill \Box$

Proposition \ref{prop:1} is one of the main results of this section. We now state some additional lemmas which will
also be useful in what follows.

\begin{lemma}
Let $\xi\in\Xi_3(V)$, $V\in\cV(n)$, and $\theta\in\Theta_{2M}$, $\theta\not\in V$.
Then $| |\xi+\theta|^{2l}-\rho^{2l}|>\rho^{2l-2+q_n}$.
\label{lem:tr}
\end{lemma}
{\em Proof.} Let $\xi=\eta+\theta'$, $\eta\in\Xi_2(V)$, $\theta'\in V\cap\dlat$.
Then $|\xi\cdot\theta|>\rho^{q_n+\epsilon_0}$ since otherwise we would have
\[
|\eta\cdot\theta|\leq |\xi\cdot\theta| +|\theta'\cdot\theta|\leq\rho^{q_{n_2}+\epsilon_0}
 +c\rho^{q_n}\leq 2\rho^{q_n+\epsilon_0}|\theta|
\]
and therefore $\eta\in\Xi_1(V+\{t\theta : t\in\R\})$ by (\ref{v1v2}). Hence
\begin{eqnarray*}
| |\xi+\theta|^{2l}-\rho^{2l}|&\geq& | |\xi+\theta|^{2l}-|\xi|^{2l}| -| |\xi|^{2l}-\rho^{2l} | \\
&\geq&c\rho^{2l-2}\rho^{q_n+\epsilon_0} -\rho^{\gamma}\\
&>&\rho^{2l-2+q_n}\, . \hspace{4cm} \Box
\end{eqnarray*}
\begin{lemma}
Let $V\in\cV(n)$ and let $\xi\in\Xi_3(V)$ and $\theta\in V\cap\dlat$. If
$\xi+\theta\not\in\Xi_3(V)$ then $| |\xi +\theta|^{2l}-\rho^{2l}|\geq\rho^{\gamma}$.
\label{nac}
\end{lemma}
{\em Proof.} Let $\xi=\eta'+\theta'$, $\eta'\in\Xi_2(V)$, $\theta'\in V\cap\dlat$.
Then $\xi+\theta=\eta' +(\theta+\theta')\in
\Xi_2(V)+(V\cap\dlat)$, and the result follows from the definition of $\Xi_3(V)$. $\hfill \Box$

\begin{corollary}
Let $\xi\in\Xi_3(V)$ and $\theta\in\Theta_{2M}$. If $\xi+\theta\not\in\Xi_3(V)$ then \nl
$| |\xi+\theta|^{2l}-\rho^{2l}|>\rho^{\gamma}$.
\label{cor:loud}
\end{corollary}
{\em Proof.} The result follows from Lemma \ref{lem:tr} if $\theta\in V$ and from Lemma \ref{nac}
if $\theta\not\in V$. $\hfill \Box$

%%%%%%%%%%  main proposition  %%%%%%%%%%%%

We can now use the results obtained so far and apply Lemma
\ref{ll2} in our context. Let $k\in\dfd$ be fixed and let the
projection $\bar{P}=\bar{P}(k)$ be as in Corollary
\ref{cor:where}. Given a lattice subspace $V\in\cV(n)$, $0\leq
n\leq d-1$, we set $P(V)=\cP^{(k)}(\Xi(V))$. 
The above statements imply that $P(V)\bar P=P(V)$.
The next proposition
provides information about the proximity of the parts of
$\sigma(H(k))$ (or $\sigma(\bar{P}H(k)\bar{P})$) and
$\sigma(\sum_V P(V)H(k)P(V))$ that lie near $\rho$. The sum is
taken over all lattice subspaces $V$.
\begin{proposition}
There exists a map $G$ from the set of all eigenvalues of the operator 
$\bar H(k):=\sum_V P(V)H(k)P(V)$ that lie in $J$ into the set
of all eigenvalues of $\bar{P}H(k)\bar{P}$ such that
whenever $\mu_l(\sum_V PH(k)P)\in J$, we have:
\[
| \mu_l(\sum_V P(V)H(k)P(V))- G(\mu_l(\sum_V P(V)H(k)P(V))) | \leq c\rho^{-2M(\gamma-\alpha)+\alpha}.
\]
This mapping is an injection and all eigenvalues of $\bar{P}H(k)\bar{P}$ inside $J_1:=[\lambda-90L,
\lambda+90L]$ have a pre-image under $G$.
\label{prop:main}
\end{proposition}
{\em Proof.} We apply Lemma \ref{ll2} to the operator $\bar{P}H(k)\bar{P}=\bar{P}H_0(k)\bar{P}+\bar{P}A(k)\bar{P}$. We use the projections
$\{P(V)\}$, where $V$ ranges over all possible lattice subspaces of dimension smaller than $d$; these are orthogonal
by Proposition \ref{prop:1}. Each $P(V)$ is further writen as a sum of orthogonal and invariant (for $H_0(k)$) projections,
$P(V)=\sum_{j=0}^MP_j(V)$, where
\begin{eqnarray*}
&& P_0(V)=\cP^{(k)}(\Xi_3(V)) \; , \\
&& P_j(V)=\cP^{(k)}((\Xi_3(V)+\Theta_j) \setminus (\Xi_3(V)+\Theta_{j-1}) ) \; , \;\; j=1,\ldots, M \, .
\end{eqnarray*}
It follows from Lemma \ref{lem:pvp} that
\begin{eqnarray*}
&& P(V_1)A(k)P(V_2)=0 \mbox{ if } V_1\neq V_2 \; , \\[0.2cm]
&& P_j(V)A(k)P_l(V)=0 \mbox{ if } |j-l|>1 \; , \\ [0.2cm]
&& P_j(V)A(k) (\bar{P}-\sum_VP(V))=0 \mbox{ if } j<M \, .
\end{eqnarray*}
Let us also check that the remaining two conditions of Lemma
\ref{ll2} are satisfied. We have
\begin{eqnarray*}
&& \sigma( (\bar{P}-\sum_VP(V))H_0(k)(\bar{P}-\sum_VP(V))) \\
&=&\{ |\xi|^{2l} \; : \; \{\xi\}=k \; , \;
 | |\xi|^{2l}-\rho^{2l}|<\rho^{2l}/2 \; , \;\xi\not\in\cup_V\Xi(V)\} \\
&\subset&\{ |\xi|^{2l} \; :  \;\xi\not\in\cA\},
\end{eqnarray*}
and therefore
\[
\dist( \sigma( (\bar{P}-\sum_VP(V))H_0(k)(\bar{P}-\sum_VP(V))) ,J)\geq 6\| \bar{P}A(k)\bar{P}\|
\]
We also note that for $1\leq j\leq M$,
\begin{eqnarray*}
\sigma( P_j(V)H_0(k)P_j(V))&=&\{ |\xi|^{2l} \; : \; \{\xi\}=k\; , \; \xi\in (\Xi_3(V)+\Theta_j)
\setminus (\Xi_3(V)+\Theta_{j-1})\}\\
&\subset&\{ |\xi|^{2l} \; : \;  \xi\in(\Xi_3(V)+\Theta_M)\setminus \Xi_3(V) \}\; .
\end{eqnarray*}
Corollary \ref{cor:loud} together with the fact that $\alpha<\gamma$ imply that
$a_j(V)\geq c\rho^{\gamma}$ and, in particular,
\[
a_j(V):=\dist( \sigma( P_j(V)H_0(k)P_j(V)) ,J)\geq 16\| \bar{P}A(k)\bar{P}\| \, ,\;\; j\geq 1 \, .
\]
Hence Lemma \ref{ll2} can be applied. We conclude that there exists an injection $G$ from the set of all eigenvalues of
$\sum_V P(V)H(k)P(V)$ that lie in $J$ into the set of all eigenvalues of $\bar{P}H(k)\bar{P}$ such that for any
$\mu_l(\sum_V P(V)H(k)P(V))\in J$ there holds
\begin{eqnarray*}
&&| \mu_l(\sum_V P(V)H(k)P(V)) -G(\mu_l(\sum_V P(V)H(k)P(V)) )| \\
&\leq& \max_V \left[ (6\|\bar{P}A(k)\bar{P}\|)^{2M+1}
\prod_{j=1}^M (a_j(V) -6\|\bar{P}A(k)\bar{P}\|)^{-2}\right] \\
&\leq& c\rho^{-2M(\gamma-\alpha)+\alpha};
\end{eqnarray*}
this completes the proof of the proposition. $\hfill \Box$

%%%%%%%  end of main proposition  %%%%%%%%

\

Let us briefly describe the aim of the next two sections.
Proposition \ref{prop:main} has led to the study of the operators $P^{(k)}(\Xi(V))H(k)P^{(k)}(\Xi(V))$, where $V$ is a lattice subspace and $k\in\dfd$. It will be proved that for fixed $k$ and $V$ we have a direct sum decomposition
\[
P^{(k)}(\Xi(V))H(k)P^{(k)}(\Xi(V)) =\bigoplus_{\xi}H(\xi)
\]
where the sum is taken over certain $\xi\in\Xi_2(V)$ with $\{\xi\}=k$ and $H(\xi)$ are certain operators. Rather than keeping $k$ fixed, we intend to study the spectrum of $H(\xi)$ as $\xi$ varies continuously.
For each $\xi$ we shall choose a specific eigenvalue
$\tilde{g}(\xi)$; the choice is such that $\tilde{g}(\xi)=|\xi|^{2l}$ in the unperturbed case. We then study
how $\tilde{g}(\xi)$ varies as $\xi$ varies. 
This turns out to depend on the location of $\xi$. If $\xi\in\cB$ (non-resonance region), then
$\tilde{g}(\xi)$ varies smoothly with $\xi$. If however $\xi\in\cD$ (resonance region), then we do not have this good dependence anymore.
In this case a new function $g(\xi)$ is introduced, which is very close to $\tilde{g}(\xi)$ and is
smooth along one direction only.

%%%%%%%%%%%%%%%%%%%%%%%%%%%%%%%%%%%%%%%%%%%%%%%%%%%%%%%%%%%%%%%%%%%%%%%%%%

%%%%%%%%%%%%%%%%%   NON RESONANCE REGION  %%%%%%%%%%%%%%%%%%%%%%%%%%%%%%%

%%%%%%%%%%%%%%%%%%%%%%%%%%%%%%%%%%%%%%%%%%%%%%%%%%%%%%%%%%%%%%%%%%%%%%%%%%

\section{Non-resonance region}

\begin{lemma}
If $\xi\in\cB$ and $\theta\in\Theta_M'$ then $| |\xi+\theta|^{2l}-\rho^{2l} |\gg\rho^{2l-2+q_1}$.
\label{lem:pr}
\end{lemma}
{\em Proof.} We have $|\xi\cdot\theta|\geq\rho^{q_1}$, since otherwise $\xi\in\Xi_1(\{t\theta:t\in\R\})$.
Hence, since $\alpha<2l-2+q_1$,
\begin{eqnarray*}
| |\xi+\theta|^{2l}-\rho^{2l}|&\geq& | |\xi+\theta|^{2l}-|\xi|^{2l}| - | |\xi|^{2l}-\rho^{2l} |\\
&\geq&c\rho^{2l-2+q_1} -100L\rho^{\alpha} \\
&\geq&c\rho^{2l-2+q_1},
\end{eqnarray*}
as required. $\hfill \Box$

Let us fix a $\xi\in\cB$ with $\{\xi\}=:k$. From (\ref{kosm}) we have
\begin{equation}
\| \cP^{(k)}(\xi+\Theta_M)A(k)\cP^{(k)}(\xi+\Theta_M)\|\leq L\rho^{\alpha}.
\label{ras}
\end{equation}
\begin{lemma}
There exists a unique eigenvalue $\tilde{g}(\xi)$ of $\cP^{(k)}(\xi+\Theta_M)H(k)\cP^{(k)}(\xi+\Theta_M)$ which lies
within distance $L\rho^{\alpha}$ of $|\xi|^{2l}$.
\label{11}
\end{lemma}
{\em Proof.} The existence of such an eigenvalue follows from (\ref{ras}) and the min-max principle.
To prove its uniqueness we argue by contradiction: let us assume that there exist two such eigenvalues.
Then there exists an eigenvalue
of $\cP^{(k)}(\xi+\Theta_M)H_0(k)\cP^{(k)}(\xi+\Theta_M)$ which is different from $|\xi|^{2l}$ and which is within
distance $2L\rho^{\alpha}$ of $|\xi|^{2l}$. Hence $| |\xi+\theta|^{2l} -|\xi|^{2l}|<2L\rho^{\alpha}$
for some $\theta\in\Theta_M'$. This implies that
\[
| |\xi+\theta|^{2l}-\rho^{2l}|\leq (102 L )\rho^{\alpha}\; ,
\]
contradicting Lemma \ref{lem:pr}. $\hfill \Box$

\

We shall obtain some more information on the eigenvalues $\tilde{g}(\xi)$, $\xi\in\cB$. First, we observe that the matrix of
$\cP^{(k)}(\xi+\Theta_M)H(k)\cP^{(k)}(\xi+\Theta_M)$ (with respect to the basis $\{e_{\xi+\theta}\}_{\theta\in\Theta_M}$
of $\ran \cP^{(k)}(\xi+\Theta_M)$ and
for some ordering $0=\theta_0,\theta_1,\theta_2,\ldots$ of $\Theta_M$) has the form
\begin{equation}
\left(
\begin{array}{cccc}
a_{00}(\xi) & a_{0\theta_1}(\xi) & a_{0\theta_2}(\xi) & \ldots \\
a_{\theta_10}(\xi) & a_{\theta_1\theta_1}(\xi) & a_{\theta_1\theta_2}(\xi)  & \ldots \\
a_{\theta_20}(\xi) & a_{\theta_2\theta_1}(\xi) & a_{\theta_2\theta_2}(\xi)  & \ldots \\
\vdots & \vdots & \vdots  & \ddots
\end{array}
\right)
\label{for}
\end{equation}
where (cf. (\ref{hk}))
\[
a_{\theta\theta'}(\xi) =\darr{|\xi+\theta|^{2l} +(2\pi)^{-d/2}\hat{a}(0,\xi+\theta),}{\theta=\theta',}
{(2\pi)^{-d/2}\hat{a}(\theta-\theta',\xi+\theta'),}{\theta\neq\theta'. }
\]
The size of this matrix is fixed and does not depend on $\rho$.
Expanding the determinant we find that the characteristic polynomial $p(\mu)$ can be written as
\begin{equation}
p(\mu)=\Big(\prod_{\theta\in\Theta_M'}( a_{\theta\theta}(\xi)-\mu)\Big)
\Big( a_{00}(\xi) -\mu +I(\xi,\mu)\Big).
\label{and}
\end{equation}
The function $I(\xi,\mu)$ is a (finite) sum
$I(\xi,\mu)=I_1+I_2+\tilde{I}_2+\ldots$ where each $I_n$ is a linear combination of terms of the form
\begin{equation}
T_n(\xi,\mu)=\frac{P_{n+1}(\xi;\theta_1,\theta_2,\ldots)}{(a_{\theta_{i_1}\theta_{i_1}}(\xi)-\mu)(a_{\theta_{i_2}\theta_{i_2}}(\xi)-\mu)\ldots
(a_{\theta_{i_n}\theta_{i_n}}(\xi)-\mu)}
\label{t}
\end{equation}
and $\tilde{I}_n$ is a linear combination of terms of the form
\begin{equation}
\tilde{T}_n(\xi,\mu)=
\frac{(a_{00}(\xi)-\mu)P_n(\xi;\theta_1,\theta_2,\ldots)}{(a_{\theta_{i_1}\theta_{i_1}}(\xi)-\mu)(a_{\theta_{i_2}\theta_{i_2}}(\xi)-\mu)
\ldots (a_{\theta_{i_n}\theta_{i_n}}(\xi)-\mu)}\;  ;
\label{tt}
\end{equation}
here $P_k$ stands for a polynomial of degree $k$ in the off-diagonal terms $a_{\theta\theta'}(\xi)$,
$\theta\neq\theta'$, of the above matrix.
We restrict our attention to $\mu$ inside the interval
$J_{\xi}:=[ |\xi|^{2l}-L\rho^{\alpha} , |\xi|^{2l}+L\rho^{\alpha}]$ where we already know that the equation
$p(\mu)=0$ has $\tilde{g}(\xi)$ as its unique solution.

\begin{lemma}
For $\theta\in \Theta_M'$ and $\mu\in J_{\xi}$ we have 
$|a_{\theta\theta}(\xi)-\mu|\gg\rho^{2l-2+q_1}$.
\label{lem:denom}
\end{lemma}
{\em Proof.} From Lemma \ref{lem:pr} we have
\begin{eqnarray*}
|a_{\theta\theta}(\xi)-\mu| &\geq& | |\xi+\theta|^{2l} -\rho^{2l} |-|\rho^{2l}-|\xi|^{2l}|- | |\xi|^{2l}-\mu|
- (2\pi)^{-d/2}|\hat{a}(0,\xi+\theta)| \\
&\geq&   c\rho^{2l-2+q_1} -100L\rho^{\alpha}-L\rho^{\alpha}-c\rho^{\alpha}\\
&\geq&c\rho^{2l-2+q_1} \; .
\end{eqnarray*}
\vskip-2em\hfill $\Box$

\

It follows from (\ref{and}) and Lemma \ref{lem:denom} that $\tilde{g}(\xi)$ is the (unique in $J_{\xi}$) solution of
\begin{equation}
a_{00}(\xi)-\mu +I(\xi,\mu)=0 \; .
\label{111}
\end{equation}

\begin{lemma}
We have
$|\partial I/\partial \mu| \ll\rho^{-(2l-2-\alpha+q_1)}$ uniformly over all $\mu\in J_{\xi}$.
\label{lem:i}
\end{lemma}
{\em Proof.} Let $T_n$ and $\tilde{T}_n$ be as in (\ref{t}) and (\ref{tt}) respectively. Using Lemma \ref{lem:denom} we obtain by a
direct computation that
\[
\Big|\frac{\partial T_n}{\partial\mu}\Big| \leq c\rho^{-(n+1)(2l-2-\alpha+q_1)} \quad , \quad
\Big|\frac{\partial\tilde{T}_n}{\partial\mu}\Big| \leq c\rho^{-(n+1)(2l-2-\alpha+q_1)}\, ;
\]
the result follows. $\hfill \Box$

\begin{proposition}
We have
\begin{equation}
\tilde{g}(\xi)=|\xi|^{2l} +G(\xi) \; , \;\; \xi\in\cB \, ,
\label{g}
\end{equation}
where $G$ is a differentiable function satisfying
\begin{eqnarray*}
\ia && |G(\xi)|\ll\rho^{\alpha} \, ;  \\
\ib && |\nabla G(\xi) |\ll\rho^{\alpha-1}.
\end{eqnarray*}
\label{lem:illi}
\end{proposition}
{\em Proof.} We shall only prove (ii), the proof of part (i) being similar and simpler.
Let us define $G$ by (\ref{g}). From (\ref{111}) we have
\begin{equation}
I(\xi , |\xi|^{2l}+G(\xi)) -G(\xi)+ (2\pi)^{-d/2}\hat{a}(0,\xi) =0 \, .
\label{fen}
\end{equation}
Defining
\[
F(\xi, t)=I(\xi, |\xi|^{2l}+t) -t + (2\pi)^{-d/2}\hat{a}(0,\xi) \; , \;\quad\; \xi\in \cB \; , \; |t|<L\rho^{\alpha} \; ,
\]
we thus obtain $F(\xi, G(\xi))=0$ on $\cB$.
From Lemma \ref{lem:i} we have $|\partial F/\partial t| \geq 1/2$, so an application of the implicit
function theorem yields that $G$ is differentiable and
\[
|\nabla G| \leq 2 \, |\nabla_{\xi}F| \; .
\]
Hence it remains to estimate the partial derivatives $\partial F/\partial\xi_i$.
Note that $I_n(\xi, |\xi|^2+t)$ is a linear combination of terms of the form
\begin{equation}
T_n(\xi,|\xi|^{2l}+t)=
\frac{P_{n+1}(\xi;\theta_1,\theta_2,\ldots)}{\prod_{j=1}^n\Big(|\xi+\theta_{i_j}|^{2l}-|\xi|^{2l}
+(2\pi)^{-d/2}\hat{a}(0,\xi+\theta_{i_j})-t\Big)}.
\label{ti}
\end{equation}
By Lemma \ref{lem:denom} each factor in the denominator is larger in absolute value than $c\rho^{2l-2+q_1}$.
Also, the derivative of each such factor with respect to $\xi_i$ does not exceed $c\rho^{2l-2}$.
Similarly we have $|P_{n+1}|\leq c\rho^{(n+1)\alpha}$ and $|\partial P_{n+1}/\partial\xi_i| \leq c\rho^{(n+1)\alpha-1}$.
These facts imply that the partial derivatives with respect to $\xi_i$ of the RHS of (\ref{ti}) are estimated by
$c\rho^{-(2l-2+q_1-\alpha)n+\alpha-q_1}$.
The argument is similar for $\tilde{I}_n(\xi, |\xi|^2+t)$ which is a linear combination of terms
\[
\tilde{T}_n(\xi,|\xi|^{2l}+t)=\frac{((2\pi)^{-d/2}\hat{a}(0,\xi)-t)
P_{n}(\xi;\theta_1,\theta_2,\ldots)}{\prod_{j=1}^n\Big(|\xi+\theta_{i_j}|^{2l}-|\xi|^{2l}+\hat{a}(0,\xi+\theta_{i_j})-t\Big)}\; .
\]
Similar calculations show that the partial derivatives with respact to $\xi$ of this expression are also smaller than
$c\rho^{-(2l-2+q_1-\alpha)n+\alpha-q_1}$. The worst estimate corresponds to $n=1$; recalling (\ref{q1}) completes the proof of (ii).
$\hfill \Box$

%%%%%%%%%%%%%%%%%%%%%%%%%%%%%%%%%%%%%%%%%%%%%%%%%%%%%%%%%%%%%%%%%%%%%%%%%%

%%%%%%%%%%%%%%%%%%%%%    RESONANCE REGION  %%%%%%%%%%%%%%%%%%%%%%%%%%%%%%%

%%%%%%%%%%%%%%%%%%%%%%%%%%%%%%%%%%%%%%%%%%%%%%%%%%%%%%%%%%%%%%%%%%%%%%%%%%

\section{Resonance region}

We shall now study the eigenvalues of $P(V)H(k)P(V)$, where $V\in\cV(n)$, $1\leq n\leq d-1$, is fixed.
Let $\xi\in\Xi_2(V)$ be given and let $k=\{\xi\}$. We define
\[
\begin{array}{ll}
\tilde{Y}(\xi)= (\xi +(V\cap\dlat))\cap \Xi_3(V)\; , & \; Y(\xi)=\tilde{Y}(\xi)+\Theta_M \; , \\[0.2cm]
P(\xi)= \cP^{(k)}(Y(\xi)) \; , & \; H(\xi)=P(\xi)H(k)P(\xi) \; ,\\[0.2cm]
H_0(\xi)=P(\xi)H_0(k)P(\xi)\; , & \; A(\xi)=P(\xi)A(k)P(\xi)\; .
\end{array}
\]
Recalling the decomposition
$\xi=\xi_V +\xi_V^{\perp}$, $\xi_V\in V$, $\xi_V^{\perp}\in V^{\perp}$, we also define
\[
r(\xi)=|\xi_V^{\perp}| \quad , \quad \xi_V' = \xi_V^{\perp} /r(\xi)\; .
\]
We note that $r(\xi)\asymp\rho$ by (\ref{paok}).
The triple $(r(\xi),\xi_V',\xi_V)$ can be thought of as cylindrical coordinates of the point $\xi\in\Xi_2(V)$.
\begin{lemma}
$\ia$ The sets $Y(\xi)$, $\xi\in\Xi_2(V)$, either coincide or are disjoint. 

$\ib$
If $Y(\xi_1)=Y(\xi_2)$ then $\xi_1-\xi_2\in V$ and, in particular, $r(\xi_1)=r(\xi_2)$.
\label{lem:r}
\end{lemma}
{\em Proof.} Assume that $Y(\xi_1)\cap Y(\xi_2)\neq\emptyset$. Then there exist $\xi_j\in\Xi_3(V)$
and $\theta_j\in\Theta_M$, $j=1,2$, such that $\xi_1+\theta_1=\xi_2+\theta_2$.
We claim that the difference $\theta_1-\theta_2$ lies in $V$. Indeed, suppose it does not. Then the relation
$\xi_2=\xi_1 +(\theta_1-\theta_2)$ together with Lemma \ref{lem:tr} yields
$||\xi_2|^{2l}-\rho^{2l}|>\rho^{2l-2+q_n}\geq\rho^{\gamma}$, contradicting the fact that $\xi_2\in\Xi_3(V)$.
Hence $\xi_2-\xi_1\in V$, and both $\ia$ and $\ib$ follow. $\hfill \Box$

\

Part $\ia$ of Lemma \ref{lem:r} points to an equivalence relation defined
on $\Xi_2(V)$, whereby $\xi_1\sim \xi_2$ if and only if $Y(\xi_1)=Y(\xi_2)$. We thus have
for each $k\in\dfd$ the direct sum decomposition
\begin{equation}
P(V)H(k)P(V)=\bigoplus H(\xi)\, ,
\label{coffee}
\end{equation}
where the sum is taken over all equivalence classes of this relation with $\{\xi\}=k$.

Hence we intend to study the operators $H(\xi)$, $\xi\in\Xi_2(V)$.
In fact, we shall compare the eigenvalues of two such operators $H(\xi_1)$ and $H(\xi_2)$;
this will be carried out using auxiliary operators denoted by $H(\xi,U)$, where $\xi\in\Xi_2(V)$ and $U$ is a subset of $\Xi_2(V)$ containing $\xi$.
We therefore introduce some additional definitions: given $\xi_1,\xi_2\in\Xi_2(V)$ and letting $k_1=\{\xi_1\}$
we set
\[
\begin{array}{ll}
Y(\xi_1,\xi_2)=Y(\xi_1)\cup (Y(\xi_2)-\xi_2+\xi_1) \; , & \; P(\xi_1,\xi_2)=\cP^{(k_1)}(Y(\xi_1,\xi_2)) \; , \\[0.2cm]
H(\xi_1,\xi_2)=P(\xi_1,\xi_2)H(k_1)P(\xi_1,\xi_2) \; , & \; H_0(\xi_1,\xi_2)=P(\xi_1,\xi_2)H_0(k_1)P(\xi_1,\xi_2) \; ,\\[0.2cm]
A(\xi_1,\xi_2)=P(\xi_1,\xi_2)A(k_1)P(\xi_1,\xi_2) \; . & \;
\end{array}
\]
Finally, given a set $U\subset\Xi_2(V)$ containing $\xi$ we define (with $k:=\{\xi\}$)
\[
\begin{array}{ll}
Y(\xi,U)=\cup_{\xi_1\in U}Y(\xi,\xi_1) \; , & \; P(\xi,U)=\cP^{(k)}(Y(\xi,U)) \; , \\[0.2cm]
H(\xi,U)=P(\xi,U)H(k)P(\xi,U) \; , & \; H_0(\xi,U)=P(\xi,U)H_0(k)P(\xi,U) \; ,\\[0.2cm]
A(\xi,U)=P(\xi,U)A(k)P(\xi,U) \; . & \;
\end{array}
\]
We also define the isometry $F_{\xi_1,\xi_2}: \ran P(\xi_1,U) \to \ran P(\xi_2,U)$ by
\[
F_{\xi_1,\xi_2}e_{\eta}=e_{\eta+\xi_2-\xi_1}\; , \;\; \eta\in Y(\xi_1,U).
\]
\begin{lemma}
Let $\xi_1,\xi_2\in\Xi_2(V)$ satisfy $|\xi_1-\xi_2| \leq c\rho^{\alpha-2l+1}$.
Then for any $\xi\in Y(\xi_1,\xi_2)\setminus Y(\xi_1)$ we have
$| |\xi|^{2l} -\rho^{2l}|> \rho^{\gamma}$.
\label{lem:ppp}
\end{lemma}
{\em Proof.} Let $\xi'=\xi+\xi_2-\xi_1$. Then $\xi'\in Y(\xi_2)$ and therefore
$\xi'=\eta+\theta$ where $\eta\in\Xi_3(V)$, $\theta\in\Theta_M$ and $\eta-\xi_2\in V\cap\dlat$.
We distinguish two cases. \nl
$\ia$ $\theta\not\in V$.  In this case Lemma \ref{lem:tr} gives
$| |\xi'|^{2l} -\rho^{2l}| > \rho^{2l-2+q_n}$. Therefore
\[
| |\xi|^{2l} -\rho^{2l}| \geq  | |\xi'|^{2l} -\rho^{2l}| - | |\xi'|^{2l} -|\xi|^{2l}| > \rho^{2l-2+q_n} -c\rho^{2l-1}|\xi_2-\xi_1| >\rho^{\gamma}.
\]
$\ib$ $\theta\in V$. Then $\xi-\xi_1=(\eta-\xi_2)+\theta\in V\cap\dlat$.
Therefore, since $\xi_1\in\Xi_2(V)$ and $\xi\not\in Y(\xi_1)\supset\Xi_3(V)$, it follows that
$| |\xi|^{2l}-\rho^{2l} |>\rho^{\gamma}$, which completes the proof in this case. $\hfill \Box$

\begin{lemma}
Let $\xi\in\Xi_2(V)$ and let $U\subset\Xi_2(V)$ be a set of diameter smaller than $c\rho^{\alpha -2l+1}$ containing $\xi$. Then
there exists an injection $G$ from the set of all eigenvalues of $H(\xi)$
into the set of all eigenvalues of $H(\xi,U)$, such that
each eigenvalue of $H(\xi, U)$ in $ J_1:=[\rho^{2l}-98L,\rho^{2l}+98L]$ is in the range of $G$. Moreover
for any $\mu_i(H(\xi))\in J_1$ we have
\[
| \mu_i(H(\xi)) - G(\mu_i(H(\xi)))| < c\rho^{-2M(\gamma-\alpha)+\alpha},
\]
and
\[
G(\mu_i(H(\xi)))=\mu_{i+l}(H(\xi,U)),
\]
where $l=:l(\xi,U)$ is the number of points $\eta\in Y(\xi,U)\setminus Y(\xi)$ such that $|\eta|<\rho$.
\label{lem:achl11}
\end{lemma}
{\em Proof.} The proof is an application of Lemma \ref{ll2} so let us verify that all its conditions
are satisfied. The lemma is applied to the operator $H(\xi,U)$ which is the
sum of $H_0(\xi,U)$ and the perturbation $A(\xi,U)$; we note that
$\|A(\xi,U)\|\leq L\rho^{\alpha}$ by Lemma \ref{lem:norm}.
We apply the lemma with $n=0$. The projection $P^0$ is $P^0=P(\xi)$ and is
decomposed as a sum of orthogonal and invariant projections, $P^0=\sum_{j=0}^MP^0_j$, where
$P^0_0=\cP^{(k)}(\tilde{Y}(\xi))$ and
\[
P^0_j=P(\xi,U)\cP^{(k)}( (\tilde{Y}(\xi)+\Theta_j)\setminus (\tilde{Y}(\xi)+\Theta_{j-1}))P(\xi,U)
\, ,\;\;\; 1\leq j\leq M\, .
\]
The fact that
\begin{eqnarray*}
&&\sigma( (P(\xi,U)-P(\xi))H_0(k)(P(\xi,U)-P(\xi)))=\\
&&\hspace{3cm} =\{ |\eta|^{2l} \; : \;  \{\eta\}=k \; , \;   \eta\in Y(\xi,U) \setminus Y(\xi) \}
\end{eqnarray*}
together with Lemma \ref{lem:ppp} yields
\[
\dist( \sigma((P(\xi,U)-P(\xi))H_0(k)(P(\xi,U)-P(\xi))) , J )\geq \rho^{\gamma}.
\]
Similarly, Corollary \ref{cor:loud} yields
\[
a_j:=\dist( \sigma( P^0_jH_0(k_1)P^0_j) , J) \geq \rho^{\gamma} \, ,\;\; j\geq 1 \, .
\]
The above imply that Lemma \ref{ll2} can be applied. We conclude that there exists a map $G$ from the set
of all eigenvalues of $H(\xi)$ into the set of all eigenvalues of $H(\xi,U)$, such that each eigenvalue of
$H(\xi,U)$ in $J_1$ is in the range of $G$ and
for any $\mu_i(H(\xi))\in J_1$
\begin{eqnarray*}
&&| \mu_i(H(\xi)) - G(\mu_i(H(\xi)))| \\
&\leq& (6\|A(\xi,U)\|)^{2M+1}\prod_{j=1}^M
(a_j-6\|A(\xi,U)\|)^{-2}\\
&\leq&c\rho^{-2M(\gamma-\alpha)+\alpha},
\end{eqnarray*}
as required. $\hfill \Box$

{\bf Remark.} In order to apply Lemma \ref{ll2} we were forced to consider the smaller interval $J_1\subset J$.
There will be more occasions where our spectral interval shall need to be reduced. Strictly speaking, this will require
the introduction of several intervals $J_1\supset J_2\supset\ldots$. In order not to overburden our notation, we shall not make this explicit from now on and we shall always use the symbol $J$ for the (possibly slightly reduced) spectral integral in hand.

\

Let $\{\eta_1,\ldots,\eta_p\}\subset\Theta_M$
be a complete set of representatives of $\Theta_M$ modulo $V$, that is for
each $\theta\in\Theta_M$ there exist unique $\eta_j\in\{\eta_1,\ldots,\eta_p\}$ and $a\in V$ such that $\theta=\eta_j+a$.
Letting $V_j=\eta_j+V$ and
\[
\Psi_j(\xi)= (\xi+ (V_j\cap\dlat))\cap Y(\xi),
\]
it follows that for each $\xi\in\Xi_2(V)$
the sets $\Psi_j(\xi)$, $j=1,\ldots,p$, are pairwise disjoint and
\[
Y(\xi)=\bigcup_{j=1}^p \Psi_j(\xi)\, .
\]

Let $U\subset\Xi_2(V)$ be a set of diameter smaller than $c\rho^{\alpha -2l+1}$ containing $\xi$.
We shall consider the matrix elements of $H(\xi,U)$ with respect to the basis
$\{e_{\eta} \, : \, \eta\in Y(\xi,U)\}$ of $\ran\, P(\xi,U)$. So let $\eta\in Y(\xi,U)$. Then there exist
a unique $\eta_k$ and a unique $\mu\in V\cap\dlat$ such that
\begin{equation}
\eta=\xi+\mu+\eta_k=r(\xi)\xi_V'+\xi_V +\mu+\eta_k \, .
\label{doc}
\end{equation}
Using Taylor's expansion we then have
\begin{eqnarray*}
|\eta|^{2l}&=&|r(\xi)\xi_V'+\xi_V +\mu+\eta_k|^{2l}\\
&=&r(\xi)^{2l}+\sum_{s=1}^{\infty}r(\xi)^{2l-s}b_s(\xi,\eta),
\end{eqnarray*}
where the function $b_s(\xi,\eta)$ has the form (using standard multi-index notation)
$b_s(\xi,\eta)=\sum_{|\alpha|=s}c_{\alpha}P_{\alpha}(\xi_V') (\xi_V+\mu+\eta_k)^{\alpha}$;
here $P_{\alpha}$ is a polynomial of degree $|\alpha|$. Hence
\begin{equation}
H(\xi,U)=r(\xi)^{2l}I + \sum_{s=1}^{\infty}r(\xi)^{2l-s}B_s(\xi,U) +A(\xi,U)\, ,
\label{deth}
\end{equation}
where the operator $B_s(\xi,U)$ is given by $B_s(\xi,U)e_{\eta}=b_s(\xi,\eta)e_{\eta}$,
$\eta\in Y(\xi,U)$. We note that Corollary \ref{cor:1} together with the fact that $\diam(U)\ll \rho^{\alpha -2l+1}$
give
\[
|\xi_V+\mu+\eta_k|=|\eta_V+(\eta_k)_V^{\perp}|\ll\rho^{q_n}.
\]
Therefore
\begin{equation}
\|B_s(\xi,U)\| =\max_{\eta}|b_s(\xi,\eta)|\ll\rho^{q_ns}\, .
\label{bug}
\end{equation}
We also note that for $s=1$ we have
\begin{equation}
b_1(\xi,\eta)=2l \xi_V'\cdot (\xi_V +\mu+\eta_k),
\label{kke}
\end{equation}
so that $\|B_1(\xi,U)\|\ll 1$.
Concerning $A(\xi,U)$, we note that for $\xi,a\in U$ and $\eta\in Y(\xi,U)$,
\[
A(\xi,U)e_{\eta}- F_{a,\xi}A(a,U)F_{\xi,a}e_{\eta}=
(2\pi)^{-d/2}\sum_{\eta'\in Y(\xi,U)} [\hat{a}(\eta'-\eta,\eta)-\hat{a}(\eta'-\eta,\eta+a-\xi)]e_{\eta}' \; .
\]
Since $|\hat{a}(\eta'-\eta,\eta)-\hat{a}(\eta'-\eta,\eta+a-\xi)| \leq c\rho^{\alpha-1}|a-\xi|$,
we can use the argument in the proof of Lemma \ref{lem:norm} to obtain
\begin{equation}
\|A(\xi,U) - F_{a,\xi}A(a,U)F_{\xi,a}\| \leq c\rho^{\alpha-1}|a-\xi|\, .
\label{perta}
\end{equation}
We shall now use the above considerations to study how the eigenvalues of $H(\xi,U)$ change as $\xi$ varies.
Because the operators $H(\xi,U)$ act on different spaces, we shall use the unitary operators $F_{\xi,a}$
to move between $\ran P(\xi ,U)$ and $\ran P(a ,U)$.
Let us denote by $\{\lambda_j(\xi,U)\}$ the eigenvalues of $H(\xi,U)$ in increasing order and
repeated according to multiplicity. By (\ref{deth})
\[
\lambda_j(\xi,U)=r(\xi)^{2l}+\nu_j(\xi,U),
\]
where $\{\nu_j(\xi,U)\}$ are the eigenvalues of the operator
\[
D(\xi,U):=\sum_{s=1}^{\infty}r(\xi)^{2l-s}B_s(\xi,U) +A(\xi,U).
\]
We first consider how $\lambda_j(\xi,U)$ varies as $r(\xi)$ varies.
\begin{lemma}
Let $\xi\in U\subset\Xi_2(V)$ where $\diam(U)<\rho^{\alpha -2l+1}$.
Assume that for $t$ close enough to $r(\xi)$ the point
$a(t):=t\xi_V'+\xi_V$ belongs in $U$. Let $\{\lambda_j(t)\}$ be the eigenvalues of $H(a(t),U)$. Then
for $t$ such that $a(t)\in U$,
\[
\lambda_j(t)=t^{2l}+\nu_j(t)
\]
where $\nu_j(t)$ is a function satisfying
\bee\label{eq:new5}
\frac{d\nu_j(t)}{dt}= O(\rho^{2l-2+q_n})\; .
\ene
\label{lem:sim}
\end{lemma}
{\em Proof.}
Replacing $\xi$ by $a(t)$ in (\ref{deth}) yields
\begin{equation}
H(a(t),U)=t^{2l}I + \sum_{s=1}^{\infty}t^{2l-s}B_s(a(t),U) +A(a(t),U) \, .
\label{deth1}
\end{equation}
Hence $\lambda_j(t)=t^{2l}+\nu_j(t)$, where $\{\nu_j(t)\}$ are the eigenvalues of the operator \nl
$\sum t^{2l-s}B_s(a(t),U)+A(a(t),U)$. Now a simple computation shows that
\[
F_{a(t),\xi}B_s(a(t),U)F_{\xi,a(t)}= B_s(\xi,U),
\]
hence
\[
F_{a(t),\xi}\Big(\sum_{s=1}^{\infty}t^{2l-s}B_s(a(t),U)\Big)F_{\xi,a(t)}=\sum_{s=1}^{\infty}t^{2l-s}B_s(\xi,U)=:B(t).
\]
We also have for $\eta\in Y(\xi,U)$,
\[
F_{a(t),\xi}A(a(t),U)F_{\xi,a(t)}e_{\eta}=(2\pi)^{-d/2}\sum_{\eta'\in Y(\xi,U)}\hat{a}(\eta'-\eta,\eta+a(t)-\xi)e_{\eta'}
=:A(t)e_{\eta}\, .
\]
We note that both $B(t)$ and $A(t)$ act on the same space which is $t$-independent -- namely $\ran\, P(\xi,U)$.
Therefore $\nu_j(t)=\mu_j(B(t)+A(t))$. Letting $\{\phi_j(t)\}$ be an orthonormal set of eigenfunctions of
$B(t)+A(t)$ and recalling (\ref{bug}) we use standard perturbation theory to obtain
\begin{eqnarray*}
\Big|\frac{d\nu_j(t)}{dt}\Big|&=&\Big|\inprod{\frac{dB(t)}{dt}\phi_j(t)}{\phi_j(t)} +
\inprod{\frac{dA(t)}{dt}\phi_j(t)}{\phi_j(t)}\Big|\\
&\leq &\|\frac{dB(t)}{dt}\| +\|\frac{dA(t)}{dt}\| \\
&\ll&\sum_{s=1}^{\infty}|2l-s|t^{2l-s-1}\|B_s(\xi,U)\|  + \rho^{\alpha-1}\\
&\ll&\rho^{2l-2+q_n}\; ;
\end{eqnarray*}
here we also used (\ref{perta}) to estimate $\|dA(t)/dt\|$.
This completes the proof. $\hfill \Box$

\

{\bf Remark.} We can slightly improve estimate (\ref{eq:new5}) if we use (\ref{kke}); this however would not be of any use in what follows. Notice that this lemma is yet another place in our paper where the proof is more complicated than in \cite{P}. 
Indeed, in \cite{P} the mappings $F_{a(t_1),a(t_2)}$ provide unitary equivalence between
$D(a(t_1),U)$ and $D(a(t_2),U)$, whereas in our paper this is no longer the case.

We next examine the case where $r(\xi)$ is fixed.
\begin{lemma}
Let $\xi,a\in U\subset\Xi_2(V)$ be such that $r(\xi)=r(a)=:r$ and $|\xi-a|<c\rho^{\alpha -2l+1}$. Then
\begin{equation}
|\lambda_j(\xi,U)-\lambda_j(a,U)|\ll\rho^{2l-2+q_n}|\xi-a| \,.
\label{eq:ma}
\end{equation}
\vspace{-.5cm}
\label{lem:hai}
\end{lemma}
{\em Proof.} Let $\eta\in Y(\xi,U)$. Using the same notation as in (\ref{doc}) we have
\begin{eqnarray*}
&& \hspace{-2cm}|B_s(\xi,U)e_{\eta}- F_{a,\xi}B_s(a,U)F_{\xi,a}e_{\eta}|\\
&=&|b_s(\xi,\eta)-b_s(a,\eta+a-\xi)| \\
&=&\Big|\sum_{|\alpha|=s}c_{\alpha}\Big[ P_{\alpha}(\xi_V')(\xi_V+\mu+\eta_k)^{\alpha}-P_{\alpha}(a_V')
(a_V+\mu+\eta_k)^{\alpha}\Big]\Big|\\
&\ll&\rho^{q_n(s-1)}|\xi-a|\, .
\end{eqnarray*}
For $s=1$ we can do better because of (\ref{kke}): we have
\[
|b_1(\xi,\eta)-b_1(a,\eta+a-\xi)| =2l|(\xi_V'-a_V')\cdot \eta_k|=\frac{2l}{r}|(\xi_V^{\perp}-a_V^{\perp})\cdot\eta_k|\ll\frac{|\xi-a| }{\rho}\, .
\]
Therefore, using also (\ref{perta}), we obtain:
\begin{eqnarray*}
&&\hspace{-1.5cm}\| F_{\xi,a}D(\xi,U)F_{a,\xi} -D(a,U)\| \\
&\leq&\sum_{s=1}^{\infty}r^{2l-s}\| F_{\xi,a}B_s(\xi,U)F_{a,\xi} -B_s(a,U)\|
+\| F_{\xi,a}A(\xi,U)F_{a,\xi} -A(a,U)\| \\
&\ll& \Big( \rho^{2l-2} +\sum_{s=2}^{\infty}[ r^{2l-s}\rho^{q_n(s-1)} ] + \rho^{\alpha-1} \Big) |\xi-a| \\
&\ll& \rho^{2l-2+q_n}|\xi-a|\, .
\end{eqnarray*}
The result follows. $\hfill \Box$

Combining the last two lemmas we have
\begin{lemma}
Let $U\subset\Xi_2(V)$ be a set with $\diam(U)\ll\rho^{\alpha -2l+1}$. Assume that $U$ contains
a piecewise $C^1$ curve joining $\xi_1,\xi_2$, of length smaller than $c|\xi_1-\xi_2|$.
Suppose that $\mu_i(H(\xi_1,U))\in J$. Then
\begin{equation}
|\mu_i(H(\xi_1,U))-\mu_i(H(\xi_2,U))| \ll\rho^{2l-1}|\xi_1-\xi_2|.
\label{nia}
\end{equation}
Suppose now in addition that $(\xi_1)_V=(\xi_2)_V$ and $(\xi_1)_V'=(\xi_2)_V'$. Then
\begin{equation}
\mu_i(H(\xi_1,U))-\mu_{i}(H(\xi_2,U)) =\{2l\rho^{2l-1} +O(\rho^{2l-2+q_n})\}(r(\xi_1)-r(\xi_2)).
\label{nia1}
\end{equation}
\label{lem:el}
\end{lemma}
\vspace{-0.8cm}
{\em Proof.} Suppose first that $(\xi_1)_V=(\xi_2)_V$ and $(\xi_1)_V'=(\xi_2)_V'$. Then by Lemma \ref{lem:sim}
there exists $t$ between $r(\xi_1)$ and $r(\xi_2)$ such that
\[
\mu_{i}(H(\xi_1,U)) -\mu_{i}(H(\xi_2,U)) =[2lt^{2l-1} +O(t^{2l-2+q_n})](r(\xi_1)-r(\xi_2)).
\]
Since $t=\rho+O(\rho^{2q_n-1})$ (cf. (\ref{paok})) estimate (\ref{nia1}) follows. Suppose next
that $r(\xi_1)=r(\xi_2)$. From Lemma \ref{lem:hai} we obtain
\[
|\mu_{i}(H(\xi_1,U)) -\mu_{i}(H(\xi_2,U))| \leq c\rho^{2l-2+q_n}|\xi_1-\xi_2|.
\]
Combining these two cases we obtain (\ref{nia}). $\hfill \Box$

\

We now proceed with some more definitions. Let $\xi\in\Xi_2(V)$ be given and $k:=\{\xi\}$.
We label the elements of $\sigma(H_0(\xi))=\{ |\eta|^{2l} \; : \; \eta\in Y(\xi)\}$ in increasing order; if there are two different points
$\eta_1,\eta_2\in Y(\xi)$ with $|\eta_1|=|\eta_2|$, then we order them in the lexicographic order of their coordinates. Hence to each $\eta\in Y(\xi)$
we have associated a natural number $j(\eta)$ such that
\[
|\eta|^{2l} =\mu_{j(\eta)}(H_0(\xi)) \; , \;\; \eta\in Y(\xi).
\]
We then define $\tilde{g}(\xi)=\mu_{j(\xi)}(H(\xi))$. It follows from Lemma \ref{lem:norm} that
\begin{equation}
|\tilde{g}(\xi)-|\xi|^{2l}|\leq L\rho^{\alpha}.
\label{synt}
\end{equation}
Let us next define for $\xi\in\Xi_2(V)$,
\[
X(\xi)=\{ \eta\in\Xi_2(V) \; : \; \eta_V=\xi_V \, , \; \eta'_V=\xi'_V\}.
\]
Clearly $X(\xi)$ is a union of at most finitely many intervals; without any loss of generality we assume that $X(\xi)$ itself
is an interval.  If $\eta_1,\eta_2\in X(\xi)$, then (cf. (\ref{zeta1}))
\[
|\eta_1-\eta_2| =|r(\eta_1)-r(\eta_2)|=\frac{| |\eta_1|^2-|\eta_2|^2 |}{r(\eta_1)+r(\eta_2)} \leq c\rho^{\alpha -2l+1},
\]
so $X(\xi)$ has length smaller than $c\rho^{\alpha -2l+1}$.

We label the elements of $\sigma(H_0(\xi,X(\xi)))=\{ |\eta|^{2l} \; : \; \eta\in Y(\xi,X(\xi))\}$ in the same way as above.
Hence to each $\eta\in Y(\xi,X(\xi))$ is associated an integer $i(\eta)$ such that
\[
|\eta|^{2l} =\mu_{i(\eta)}(H_0(\xi,X(\xi))) \; , \;\; \eta\in Y(\xi,X(\xi)).
\]
We then define $g(\xi)=\mu_{i(\xi)}(H(\xi,X(\xi)))$.
Clearly $|g(\xi)-|\xi|^{2l}|\leq L\rho^{\alpha}$ for all $\xi\in\Xi_2(V)$.

\begin{lemma}
For each $\xi\in\Xi_2(V)$ the function $i(\cdot)$ is constant on $X(\xi)$.
\label{lem:ifixed}
\end{lemma}
{\em Proof.} Let $\xi_1\in X(\xi)$. We must show that the number of points of the set
$\{\eta\in Y(\xi, X(\xi))  : \, |\eta|<|\xi|\}$ coincides with the
number of points of the set $\{\eta_1\in Y(\xi_1, X(\xi)) : \, |\eta_1|<|\xi_1|\}$.
Let $\eta\in Y(\xi, X(\xi))$ be given and define $\eta_1=\eta+\xi_1-\xi$; then $\eta_1\in Y(\xi_1, X(\xi))$.
We claim that $|\eta|<|\xi|$ if and only if $|\eta_1|<|\xi_1|$. To prove this we distinguish two cases:\nl
$\ia$ $\xi-\eta\in V$. In this case $\xi_1-\eta_1\in V$, therefore
\[
|\xi_1|^2-|\eta_1|^2=|(\xi_1)_V|^2-|(\eta_1)_V|^2 =|\xi_V|^2-|\eta_V|^2=|\xi|^2-|\eta|^2 ,
\]
and the claim follows.\nl
$\ib$ $\xi-\eta\not\in V$. We shall prove that in this case
\begin{equation}
||\eta|^{2l}-\rho^{2l}|>\rho^{\gamma}.
\label{app}
\end{equation}
If $\eta\not\in Y(\xi)$ then (\ref{app}) follows from Lemma \ref{lem:ppp}, so let us assume that
$\eta\in Y(\xi)$. We then have $\eta=\bar{\eta}+\theta$ for some $\bar{\eta}\in \Xi_3(V)$ with
$\bar{\eta}-\xi\in V\cap\dlat$ and some $\theta\in\Theta_M$. Then $\theta\not\in V$ and (\ref{app}) follows from
Lemma \ref{lem:tr}. Similarly we have $||\eta_1|^{2l}-\rho^{2l}|>\rho^{\gamma}$.
Suppose now that $|\xi|<|\eta|$. Then $|\eta|^{2l}>\rho^{2l}+\rho^{\gamma}$. Hence we have
\[
|\eta_1|^{2l}\geq |\eta|^{2l} - c\rho^{2l-1}|\eta-\eta_1|\geq|\eta|^{2l}>\rho^{2l}+\rho^{\gamma} -c\rho^{\alpha}
\]
and therefore, since $\xi_1\in\cA$, we conclude that $|\eta_1|>|\xi_1|$. This completes the proof. $\hfill \Box$

\begin{lemma}
Let $\xi\in\Xi_2(V)$. Then:
\begin{eqnarray*}
\ia &&  | g(\xi)-\tilde{g}(\xi)|\ll\rho^{-2M(\gamma-\alpha)+\alpha} ; \\
\ib && g(\xi)=r(\xi)^{2l}+ s(\xi) \mbox{ where $s(\xi)$ is differentiable with respect to
$r=r(\xi)$ and}\\
&&\parder{s}{r}=O(\rho^{2l-2+q_n}) \; .
\end{eqnarray*}
\label{lem:g}
\end{lemma}
\vspace{-.8cm}
{\em Proof.} $\ia$ We apply Lemma \ref{lem:achl11} with $U=X(\xi)$. We
conclude that that there exists an injection $G$ from the set of all eigenvalues of $H(\xi)$
into the set of eigenvalues of $H(\xi,X(\xi))$ such that each eigenvalue of $H(\xi,X(\xi))$
inside $J$ belongs in the range of $G$ and for each $\mu_i(H(\xi))\in J$ we have
$|G(\mu_i(H(\xi)))-\mu_i(H(\xi))| <c\rho^{-2M(\gamma-\alpha)+\alpha}$
and moreover
\begin{equation}
G(\mu_i(H(\xi)))=\mu_{i+m}(H(\xi,X(\xi)))
\label{ear}
\end{equation}
where $m$ is the number of eigenvalues of $[P(\xi,X(\xi)) -P(\xi)]H_0(k)[P(\xi,X(\xi)) -P(\xi)]$
that are smaller than $\rho^{2l}$.
Now, it follows from the above definitions that the difference $i(\xi)-j(\xi)$ is equal to the number of points $\eta\in Y(\xi,X(\xi))\setminus Y(\xi)$
such that $|\eta|\leq |\xi|$. Because of Lemma \ref{lem:ppp}, this can be rephrased as
\begin{equation}
i(\xi)-j(\xi) = \# \{ \eta\in Y(\xi,X(\xi))\setminus Y(\xi) \;\; : \;\;  |\eta|\leq \rho \}=m\, .
\label{piano}
\end{equation}
Choosing $i=j(\xi)$ in (\ref{ear}) proves $\ia$. \nl
$\ib$ This is a direct consequence of Lemma \ref{lem:sim} (applied for $U=X(\xi)$) and Lemma \ref{lem:ifixed}. $\hfill \Box$

The fact that $g$ does not exhibit good behaviour in $\cD$ except in (locally) one direction, prevents us from
estimating $|g(b)-g(a)|$ in terms of $|b-a|$. The next lemma compensates for this; it establishes the existence of
a conjugate point $b+n$, $n\in\dlat$, which can be used in the place of $b$.
\begin{lemma}
Let $[a,b]\subset\Xi_2(V)$ be a segment of length $|b-a|<c\rho^{\alpha -2l+1}$. Then there exists
$n\in\dlat$ such that
$|g(b+n)-g(a)|\leq c\rho^{2l-1}|b-a| +O(\rho^{-2M(\gamma-\alpha)+\alpha})$.
Suppose now in addition that there exists $m\in\dlat\setminus\{0\}$ such that
$[a+m,b+m]\subset\Xi_2(V)$. Then there exists $n_1\in\dlat$, $n_1\neq n$, such that
$|g(b+n_1)-g(a+m)|\leq c\rho^{2l-1}|b-a| +O(\rho^{-2M(\gamma-\alpha)+\alpha})$.
\label{lem:ath}
\end{lemma}
{\em Proof.} The proof is almost identical to the proof of lemma 7.11 from \cite{P},
so we will skip it. $\hfill \Box$

We can now state the following lemma, which collects together the previous results.

\begin{lemma}
Let $V\in\cV(n)$, $1\leq n\leq d-1$, and $M>0$ be given. There exist mappings $g,\tilde{g}:\Xi_2(V)\to\R$ with the following properties:
\begin{eqnarray*}
&\ia& \mbox{$\tilde{g}(\xi)$ is an eigenvalue of $P(V)H(k)P(V)$, where $k:=\{\xi\}$. Moreover, for}\\
&& \mbox{each $k$, all eigenvalues of  $P(V)H(k)P(V)$ inside $J$ are in the image of $\tilde{g}$;} \\
&\ib& \mbox{If $\xi\in\cA$, then}\\
&& \;\; {\rm (a)}\quad |\tilde{g}(\xi)-g(\xi)|\ll\rho^{-2M(\gamma-\alpha)+\alpha}; \\
&& \;\; {\rm (b)}\quad |g(\xi)-|\xi|^{2l}|\leq 2L\rho^{\alpha} ; \\
&\ic& \mbox{$g(\xi)=r(\xi)^{2l}+ s(\xi)$ where $s(\xi)$ is differentiable with respect to
$r=r(\xi)$ and}\\
&&\parder{s}{r}=O(\rho^{2l-2+q_n}) \; .
\end{eqnarray*}
\label{lem:move}
\end{lemma}
\vspace{-0.3cm}
{\em Proof.} The first statement of (i) follows immediately from the definition of $\tilde{g}$ and (\ref{coffee}).
Parts (ii)(a)  and (iii) are contained in Lemma \ref{lem:g}.
Finally (ii)(b) follows from (\ref{synt}) and (ii)(a). $\hfill \Box$

\

We now proceed to combine the results obtained so far in this section with those
of Section 4. For this we shall need to extend the definition of $g(\xi)$
for $\xi\in\Xi_2(\{0\})=\cB$. We recall the $\tilde{g}(\xi)$ has already been defined for such $\xi$ (cf. Lemma \ref{11}). We extend $g$ in $\cB$
defining
\begin{equation}
g(\xi)=\tilde{g}(\xi)\;\; , \quad \xi\in\cB \, .
\label{kara}
\end{equation}
Hence $g$ is now a function defined on the whole of the spherical layer $\cA$.

We shall define one more function $f$ on $\cA$; this will take values in $\sigma(H)$. Let $\xi\in\cA$ be given and $\{\xi\}=:k$.
Then there exists a unique lattice subspace $V$ containing $\xi$ (so $V\in\cV(n)$ for some $n\in\{0,1,\ldots,d-1\}$; if $\xi\in\cB$ then $n=0$,
while if $\xi\in\cD$ then $n\geq 1$). As we have seen $\tilde{g}(\xi)$ is an eigenvalue of $H(\xi)$; hence (cf. (\ref{coffee})) it is an eigenvalue of
$\sum_VP(V)H(k)P(V) +QH(k)Q$.
Ordering the eigenvalues of $\sum_VP(V)H(k)P(V) +QH(k)Q$
in the usual way determines a number $\tau(\xi)\in\N$ such that $\tilde{g}(\xi)=\mu_{\tau(\xi)}(\sum_VP(V)H(k)P(V) +QH(k)Q)$. We then define
\[
f(\xi)=\mu_{\tau(\xi)}(H(k))\; .
\]
We have the following
\begin{proposition}
Let $N>0$ be given. There exist two mappings $f,g:\cA\to\R$ with the following properties:
\begin{eqnarray*}
&\ia& \mbox{$f(\xi)$ is an eigenvalue of $H(k)$, where $k:=\{\xi\}$;}\\
&\ib& \mbox{For any $k$, all eigenvalues of $H(k)$ inside $J$ are in the range of $f$}; \\
&\ic& \mbox{If $\xi\in\cA$ then $|f(\xi)-g(\xi)|\leq \rho^{-N}$ }; \\
&\id& |f(\xi)-|\xi|^{2l}| \leq c\rho^{\alpha} \, ; \\
&{({\rm v})}& \mbox{Considering the disjoint union $\cA=\cB \cup \bigcup_{n=1}^{d-1}\bigcup_{V\in\cV(n)}\Xi_2(V)$ we have:}\\
&& \!\!\!\mbox{\rm{(a)} If $\xi\in\cB$ then $g(\xi)\! =|\xi|^{2l}\! +G(\xi)$, where } |\nabla G(\xi)|\leq c\rho^{\alpha-1}\, ; \\
&&\!\!\!\mbox{\rm{(b)} If $\xi\in\Xi_2(V)$ then $g(\xi)=r(\xi)^{2l}+s(\xi)$, where } \Big|\parder{s}{r}\Big|\leq c\rho^{2l-2+q_n}.
\end{eqnarray*}
\label{prop:move}
\end{proposition}
{\em Note.} We do not claim -- and indeed it is not the case in general
-- that either $f$ or $g$ is continuous in $\cA$. \newline
{\em Proof.} Part (i) is trivial. Part (ii) follows from Lemma \ref{ll2}. 
The same lemma together with Corollary
\ref{cor:where} implies that
\begin{equation}
|f(\xi)-\tilde{g}(\xi)|\ll\rho^{-2M(\gamma-\alpha)+\alpha}.
\label{pli}
\end{equation}
This, together with Lemma \ref{lem:move} (ii)(a) (if $\xi\in\cD$) or (\ref{kara}) (if $\xi\in\cB$), implies (iii) if we choose $M$ sufficiently large so that 
$-2M(\gamma-\alpha)+\alpha<-N$.
Part (iv) follows from (iii) and Lemma \ref{lem:move} (ii). Finally parts (v) (a) and (b) follow from Proposition \ref{lem:illi}
and Lemma \ref{lem:g} (ii) respectively. $\hfill \Box$

The next lemma is a global version of Lemma \ref{lem:ath}; once again, the proof is
almost identical to the proof of Lemma 7.14 from \cite{P}, so we will skip it.
\begin{lemma}
Let $[a,b]\subset\cA$ be a segment of length $|b-a|<c\rho^{\alpha -2l+1}$. Then there exists $n\in\dlat$ such that
$|g(b+n)-g(a)|\leq c\rho^{2l-1}|b-a| +O(\rho^{-2M(\gamma-\alpha)+\alpha})$.
Suppose now in addition that there exists $m\in\dlat\setminus\{0\}$ such that
$[a+m,b+m]\subset\cA$. Then there exists $n_1\in\dlat$, $n_1\neq n$, such that
$|g(b+n_1)-g(a+m)|\leq c\rho^{2l-1}|b-a| +O(\rho^{-2M(\gamma-\alpha)+\alpha+d})$.
\label{lem:athglob}
\end{lemma}

%%%%%%%%%%%%%%%  FINAL SECTION   %%%%%%%%%%%%%%

\section{Proof of the main theorem}

Let $\delta$, $0<\delta\le\rho^{2l-3}$, be a parameter, the precise value of which will be determined later on.  
We denote by $\cA(\delta)$, $\cB(\delta)$ and $\cD(\delta)$ the intersections of
$g^{-1}( [\rho^{2l}-\delta ,\rho^{2l}+\delta])$ with $\cA$,
$\cB$ and $\cD$ respectively. 
\begin{lemma}
There holds
\begin{eqnarray*}
&\ia&\vol(\cA(\delta))\asymp\delta\rho^{d -2l} \; ,\\
&\ib&\vol(\cB(\delta))\asymp\delta\rho^{d -2l} \; ,\\
&\ic&\vol(\cD(\delta))\leq\delta\rho^{d-2l-\epsilon_0} \, ,
\end{eqnarray*}
provided $\rho$ is large enough.
\label{lem:volumes}
\end{lemma}
{\em Proof.} It is enough to prove (ii) and (iii). Let us consider a point $\xi\in\cB$.
We write $\xi=r\xi'$ where $r>0$ and $|\xi'|=1$. Definition (\ref{kara}) together with Proposition \ref{lem:illi} implies that
\[
\parder{g}{r} \asymp \rho^{2l-1},
\]
uniformly over $\xi\in\cB$. Hence for each $\xi'$ the segment
\begin{equation}
\{ \xi= r\xi'\in\cB \; : \; r>0 \; , \; g(\xi)\in  [\rho^{2l}-\delta ,\rho^{2l}+\delta] \}
\label{xt}
\end{equation}
is an interval of length $\asymp\delta\rho^{-2l+1}$. Integration over all $\xi'\in S^{d-1}$ yields (ii).
To prove (iii), let us consider a point $\xi\in\Xi_2(V)$ and let $(r,\xi_V',\xi_V)$ be the
corresponding cylindrical coordinates. 
For $\theta\in\Theta_{6M}'$ let
\[
\cD_{\theta}(\delta)=\{\xi\in\cA(\delta) \; : \;  |\xi\cdot\theta |\leq \rho^{1-\epsilon_0}|\theta| \}.
\]
It follows from (v) of Proposition \ref{prop:move} that
\[
\parder{g}{r} \asymp \rho^{2l-1}.
\]
Thus, the intersection of $\cD_{\theta}(\delta)$ with 
the semi-infinite interval  $\{\xi=(r,\xi_V',\xi_V), \, r>0\}$, with $(\xi_V',\xi_V)$
being fixed, is an interval of length smaller than $c\delta\rho^{-2l+1}$.
Therefore,
$\vol(\cD_{\theta}(\delta))\ll (\delta\rho^{-2l+1})\rho^{1-\epsilon_0} \rho^{d-2}$.
The number of points $\theta\in\Theta_{6M}'$ is fixed. Hence, since $\cD(\delta)\subset\cup_{\theta\in\Theta_{6M}'}
\cD_{\theta}(\delta)$, we conclude that
\[
\vol(\cD(\delta))\leq\sum_{\theta\in\Theta_{6M}'}\vol(\cD_{\theta}(\delta))\ll \delta\rho^{d-2l-\epsilon_0}\, ,
\]
which implies (iii). $\hfill \Box$

The next lemma is crucial for the proof of the main theorem. It gives an upper estimate on the volume
of intersections of translates of $\cB(\delta)$. Recall that when $\bxi\in\cA$, 
we have $g(\xi)=|\xi|^{2l}+G(\xi)$, where $|G(\xi)|=O(\rho^{\al})$ and
$|\nabla G(\xi)|=O(\rho^{\al-1})$.
\begin{lemma}

(i) Let $d\geq 3$. Then
\begin{equation}
\vol\Big(\cB(\delta)\cap (\cB(\delta)+a)\Big) \ll (\delta^2\rho^{-4l+d+1} + \delta\rho^{(\alpha-2l+1)d-\alpha}) \; ,
\label{static}
\end{equation}
uniformly over all $a\in\R^d$ with $|a|\geq C$ for any positive constant $C$. In addition, there exists $c_4>0$ such that if $a$ satisfies $| |a|-2\rho|\geq c_4$, then
\begin{equation}
\vol\Big(\cB(\delta)\cap (\cB(\delta)+a)\Big) \leq c\delta^2\rho^{-4l+d+1}.
\label{static1}
\end{equation}

(ii) If $d=2$, then 
\bee\label{eq:staticnew2}
\vol\Big(\cB(\delta)\cap (\cB(\delta)+a)\Big) \ll\de^{3/2}\rho^{3-3l}+\de\rho^{\alpha -4l+2}.
\ene
In addition, there exists $c_4>0$ such that if $a$ satisfies $| |a|-2\rho|\geq c_4$, then
\bee\label{eq:staticnew3}
\vol\Big(\cB(\delta)\cap (\cB(\delta)+a)\Big) \ll\de^{3/2}\rho^{3-3l}.
\ene
\label{lem:volinters}
\end{lemma}
{\em Proof.} First of all, we notice that it is enough to prove this lemma assuming
$l=1$. Indeed, suppose we have established Lemma \ref{lem:volinters} for $l=1$. In the
general case, we introduce a new function $\tilde g(\xi):=(g(\xi))^{1/l}$. Then
a simple calculation shows that $\tilde g(\xi)=|\xi|^2+\tilde G(\xi)$, where $|\tilde G(\xi)|=O(\rho^{\tilde \al})$ and
$|\nabla \tilde G(\xi)|=O(\rho^{\tilde \al-1})$, with $\tilde\al=\al+2-2l<1$. 
Moreover, 
\bee
\cB(\delta)\subset\tilde\cB(\tilde\delta):=\{\xi\in\cB,\,\tilde g(\xi)\in [\rho^2-\tilde\delta,\rho^2+\tilde\delta]\},
\ene
with $\tilde\delta=\delta\rho^{2-2l}$. Thus, applying (\ref{static}) for $\tilde g$
(with $l=1$), we obtain:
\begin{equation}
\vol\Big(\cB(\delta)\cap (\cB(\delta)+a)\Big) \ll (\tilde\delta^2\rho^{d-3} + \tilde\delta\rho^{(\tilde\alpha-1)d-\tilde\alpha}).
\label{staticnew}
\end{equation}
After inserting the expressions defining $\tilde\delta$ and $\tilde\alpha$, we obtain
(\ref{static}). The other estimates are similar. 

Thus, from now on we assume without loss of generality that $l=1$. In this case we need to prove the following estimates:
\begin{equation}
\vol\Big(\cB(\delta)\cap (\cB(\delta)+a)\Big) \ll (\delta^2\rho^{d-3} + \delta\rho^{(\alpha-1)d-\alpha}),\,\,\, d\ge 3;
\label{staticnew1}
\ene
\bee\label{staticnew2}
\vol\Big(\cB(\delta)\cap (\cB(\delta)+a)\Big) \ll\de^{3/2}+\de\rho^{\alpha-2},\,\,\, d=2,
\ene
and if $| |a|-2\rho|\geq c_4$, then
\bee\label{staticnew3}
\vol\Big(\cB(\delta)\cap (\cB(\delta)+a)\Big) \ll \delta^2\rho^{d-3},\,\,\, d\ge 3,
\ene
\bee\label{staticnew4}
\vol\Big(\cB(\delta)\cap (\cB(\delta)+a)\Big) \ll\de^{3/2},\,\,\, d=2.
\ene
Denote by $C_2$ a constant such that 
\bee\label{G}
|G(\xi)|\le C_2\rho^{\al},\,\,\,\,|\nabla G(\xi)|\le C_2\rho^{\beta},
\ene
where we have denoted 
\bee
\beta:=\al-1<0.
\ene 
We need to estimate the volume of the set
\[
\CX:=\{\bnu\in\bigl(\CB\cap(\CB+a)\bigr)\, : \;  g(\bnu)\in [\rho^2-\de,\rho^2+\de], \;\; g(\bnu-a)\in [\rho^2-\de,\rho^2+\de]\}.
\]
First, we will estimate the $2$-dimensional area of the intersection of $\CX$ with an arbitrary
$2$-dimensional plane containing the origin and the vector $a$; the volume of $\CX$ then can
be obtained using the integration in cylindrical coordinates. So, let $\GV$ be any
$2$-dimensional plane containing the origin and $a$, and let us estimate
the area of $\CX_\GV:=\GV\cap\CX$.
Let us introduce cartesian coordinates in $\GV$ so that
$\bnu\in\GV$ has coordinates $(\nu_1,\nu_2)$
with $\nu_1$ going along $a$, and $\nu_2$ being orthogonal to $a$.
For any $\bnu\in\CX_\GV$ we have
\bees
g(\bnu)=\nu_1^2+\nu_2^2+G(\bnu),
\enes
and so
\bees
2\de\ge |g(\bnu)-g(\bnu-a)|=|\nu_1^2-(\nu_1 -|a|)^2+G(\bnu)-G(\bnu-a)|=|a|\ |(|a|-2\nu_1)|+|a|O(\rho^{\beta}).
\enes
This implies that $|\frac{|a|}{2}-\nu_1|=O(\rho^{\beta})+O(\rho^{-1})$ and, therefore,
\bee\label{wherenu1is}
\frac{|a|}{3}<\nu_1<\frac{2|a|}{3},
\ene
whenever $\bnu\in\CX_\GV$. This implies
\bee\label{partialg1}
\frac{\partial g(\bnu)}{\partial\nu_1}\gg |a| \; .
\ene
Thus, for any fixed $t\in\R$,
the intersection of the line $\nu_2=t$ with $\CX_\GV$ is an interval
of length $\ll \rho^{\beta}$ (we can assume without loss of generality that $\beta>-1/2$).

Let us cut $\CX_\GV$ into two parts: $\CX_\GV=\CX_\GV^1\cup\CX_\GV^2$ with
$\CX_\GV^1:=\{\bnu\in\CX_\GV,\,|\nu_2|\le 2C_2\rho^{\beta}\}$,
$\CX_\GV^2=\CX_\GV\setminus\CX_\GV^1$, and estimate the volumes of these sets ($C_2$ is the
constant from \eqref{G}). We start from $\CX_\GV^1$. Suppose that $\CX_\GV^1$ is non-empty,
say $\bxi=(\nu_1,\nu_2)\in\CX_\GV^1$ (note that $|\nu_2|\ll\rho^{\beta}$). 
Denote $\boldeta:=(\nu_1,0)$. Then 
\bee\label{vol1}
g(\boldeta)=g(\bxi)+O(\rho^{2\beta})=\rho^2+O(\rho^{2\beta}).
\ene
Similarly, 
\bee\label{vol2}
g(\boldeta-a)=g(\bxi-a)+O(\rho^{2\beta})=\rho^2+O(\rho^{2\beta}).
\ene
Thus, if $\CX_\GV^1$ is non-empty, then $|a|\sim\rho$ and the first coordinate of any point $\bxi\in\CX_\GV^1$ 
satisfies $\nu_1\sim\rho$. Therefore, we have $\frac{\partial g}{\partial\nu_1}(\bxi)\gg \rho$. Let $s_1$ denote
the unique positive solution of the equation $g(s_1,0)=\rho^2$; similarly, let $s_2$ be the unique positive solution 
of the equation $g(-s_2,0)=\rho^2$.
Our conditions on $g$ imply that $s_j=\rho+O(1)$. 
Estimate (\ref{vol1}) implies $\nu_1=s_1+O(\rho^{2\beta-1})$; similarly, 
estimate (\ref{vol2}) implies $\nu_1-|a|=-s_2+O(\rho^{2\beta-1})$. Thus, 
if $\CX_\GV^1$ is non-empty, then we have
\bee\label{vol3}
|a|=s_1+s_2+O(\rho^{2\beta-1})=2\rho+O(1).
\ene
Let us now fix $t$, $0\le t\le 2C_2\rho^{\beta}$. Since $\frac{\partial g}{\partial\nu_1}(\bxi)\gg \rho$,
the length of the intersection
of $\CX_\GV^1$ with the line $\nu_2=t$ is $\ll\de\rho^{-1}$. This 
implies that the area of $\CX_\GV^1$ is $\ll \rho^{\beta-1}\de$.
Now we define the `rotated' set $\CX^1$ which consists of the points from $\CX$
which belong to $\CX_\GV^1$ for some $\GV$. Computing the volume of this set using
integration in the cylindrical coordinates, we obtain
\bee\label{CX1}
\volume(\CX^1)\ll \rho^{(d-1)\beta-1}\de.
\ene
We now consider $\CX_\GV^2$. Let us decompose $\CX_\GV^2=\overline{\CX_\GV^2}\cup\underline{\CX_\GV^2}$, where
\bees
\overline{\CX_\GV^2}=\{\bnu\in\CX_\GV^2:\nu_2>0\}
\enes
and
\bees
\underline{\CX_\GV^2}=\{\bnu\in\CX_\GV^2:\nu_2<0\}.
\enes
Notice that for any $\bnu\in\overline\CX_\GV^2$, formula (\eqref{G}) implies
\bee\label{partialg2}
\frac{\partial g(\bnu)}{\partial\nu_2}\gg \nu_2.
\ene
Let $\bnu^l=(\nu^l_1,\nu^l_2)$ be the point in
the closure of $\overline{\CX_\GV^2}$
with the smallest value of the first coordinate:
$\nu^l_1\le\nu_1$ for any $\bnu=(\nu_1,\nu_2)\in\overline{\CX_\GV^2}$.
Analogously, we define $\bnu^r$
to be the point in the closure of $\overline{\CX_\GV^2}$
with the largest first coordinate,
$\bnu^t$ the point with the largest second coordinate, and $\bnu^b$ the point
with the smallest second coordinate. %(see Figure \ref{fig:9} for an illustration). 
Note that $\nu^t_2\ll\rho$.
%\begin{figure}[!hbt]
%\begin{center}
%\framebox[0.8\textwidth]{\includegraphics[width=0.60\textwidth]{lp_pict9}}
%\caption{The set $\overline{\CX_\GV^2}$ (the area bounded by four arcs)\label{fig:9}}
%\end{center}
%\end{figure}

Let us prove that
\bee\label{width0}
\nu_1^r-\nu^l_1\ll \de.
\ene

Indeed, suppose first that $\nu^r_2\ge \nu^l_2$. Let $\bnu^{rl}:=(\nu^r_1,\nu^l_2)$. Then,
since $g$ is an increasing function of $\nu_2$ when $\nu_2>2C_2\rho^{\beta}$, we have
$g(\bnu^{rl})\le g(\bnu^r)\le \rho^2+\de$. Therefore,
$g(\bnu^{rl})-g(\bnu^l)\le 2\de$. Since $g$ is increasing with respect to $\nu_1$,
estimate (\ref{partialg1}) implies \eqref{width0}.

Suppose now that $\nu^r_2\le \nu^l_2$. Let $\bnu^{lr}:=(\nu^l_1,\nu^r_2)$. Then
$g(\bnu^{lr}-a)\le g(\bnu^l-a)\le \rho^2+\de$.
Therefore,
$g(\bnu^{lr}-a)-g(\bnu^r-a)\le 2\de$. Since $g(\cdot - a)$ is decreasing with respect to $\nu_1$,
\eqref{wherenu1is} and \eqref{partialg1} imply \eqref{width0}.

Thus, we have estimated the width of $\overline{\CX_\GV^2}$. Let us estimate its height
(i.e. $\nu_2^t-\nu^b_2$). Let us assume that $\nu^t_1\ge\nu^b_1$; otherwise, we use the same
trick as in the previous paragraph and consider $g(\cdot-a)$ instead of $g$.
Let $\bnu^{bt}:=(\nu^b_1,\nu^t_2)$. Then
$g(\bnu^{bt})\le g(\bnu^t)\le \rho^2+\de$. Therefore,
$g(\bnu^{bt})-g(\bnu^b)\le 2\de$. Now, \eqref{partialg2} implies
\bee\label{hightd2}
(\nu_2^t)^2-(\nu_2^b)^2=2\int_{\nu_2^b}^{\nu_2^t}\nu_2d\nu_2
\ll \int_{\nu_2^b}^{\nu_2^t}\frac{\partial g}{\partial\nu_2}(\nu_1^b,\nu_2)d\nu_2
\le 2\de.
\ene
Therefore, we have the following estimate for the height of $\overline{\CX_\GV^2}$:
\bee\label{hight}
\nu_2^t-\nu_2^b\ll \frac{\de}{\nu_2^t+\nu_2^b}.
\ene
Now, we can estimate the volume of $\CX^2:=\CX\setminus\CX^1$ using estimates \eqref{width0}
and \eqref{hight}. Cylindrical integration produces the following:
\bee\label{CX2}
\volume(\CX^2)\ll \frac{\de^2}{\nu_2^t+\nu_2^b}(\nu_2^t)^{d-2}\le
\de^2 (\nu_2^t)^{d-3}\le \de^2\rho^{d-3}.
\ene
Equations \eqref{CX1} and \eqref{CX2} imply \eqref{staticnew1}. If $d=2$, we have
to notice that \eqref{hightd2} implies $\nu_2^t-\nu_2^b\ll \de^{1/2}$ and then use \eqref{CX1}
and \eqref{width0}. Finally, the remark before \eqref{vol3} implies 
(\ref{staticnew2}) and (\ref{staticnew4}).

$\hfill \Box$

Let $N$ be an arbitrary natural number, the precise value of which will be determined later (as well as the precise value of $\delta$). We construct the mappings $f,g$ according to Proposition \ref{prop:move}.
\begin{lemma}
Let $\xi\in\cB(\delta)$ be a point of discontinuity of $f$. Then there exists $n\in\dlat\setminus\{0\}$ such that
$\xi+n\in\cA$ 
and
\begin{equation}
| g(\xi+n)-g(\xi)| \leq 2\rho^{-N}.
\label{1234}
\end{equation}
\label{lem:disco}
\end{lemma}
{\em Proof.} Let $\xi\in\cB(\delta)$ be a point of discontinuity of $f$. Since $f$ is bounded,
there exist two sequences $(\xi_j)$ and $(\tilde{\xi}_j)$ in $\cB(\delta)$ both converging to $\xi$
and such that the limits $\lambda:=\lim f(\xi_j)$ and $\tilde{\lambda}:=\lim f(\tilde{\xi}_j)$ both exist
in $\R$ and $\lambda\neq\tilde{\lambda}$. Let $k:=\{\xi\}$, $k_j:=\{\xi_j\}$. Since $f(\xi_j)$ are eigenvalues
of $H(k_j)$, the limit $\lambda$ is an eigenvalue of $H(k)$ \cite{K}; similarly $\tilde{\lambda}$ is an eigenvalue
of $H(k)$. Since $\lambda\neq\tilde{\lambda}$ at least one of $\lambda,\tilde{\lambda}$ is different from
$f(\xi)$, say $\tilde{\lambda}\neq f(\xi)$. The fact that $\tilde{\lambda}$ is an eigenvalue
of $H(k)$ inside $J$ implies, by (ii) of Proposition \ref{prop:move}, that $\tilde{\lambda}=f(\tilde{\xi})$
for some $\tilde{\xi}\in\cA$ with $\{\tilde{\xi}\}=k$. Using the continuity of $g$ in $\cB$ and (iii) of
Proposition \ref{prop:move},
we conclude that
\[
|g(\tilde{\xi})-g(\xi)| \leq |g(\tilde{\xi})-f(\tilde{\xi})| +  \lim | f(\tilde{\xi}_j)- g(\tilde{\xi}_j)|
\leq 2\rho^{-N},
\]
which is (\ref{1234}). We have $\tilde{\xi}\neq\xi$ since otherwise we would obtain $f(\xi)=f(\tilde{\xi})=
\tilde{\lambda}$. $\hfill \Box$

\begin{lemma}
There exists a constant $c_2>0$ with the following property:
suppose that $I\subset\cB(\delta)$ is a straight segment of length $T<\delta\rho^{-2l+1}$.
Suppose also
that there exists a point $\xi_0\in I$ such that for each $n\in\dlat\setminus\{0\}$ we have
\begin{eqnarray*}
\mbox{either }&\ia& \xi_0+n \not\in\cA \;\; 
\\
\mbox{ or } &\ib& | g(\xi_0+n)-g(\xi_0)| \geq c_2(T\rho^{2l-1}+\rho^{-N}).
\end{eqnarray*}
Then the restriction $f|_I$ is continuous.
\label{lem:c3}
\end{lemma}
{\em Proof.} We argue by contradiction. So let us assume the contrary: for any $c_2>0$ there exists a segment
$I\subset\cB(\delta)$ of length $T<\delta\rho^{-2l+1}$ and a point $\xi_0\in I$ such that for any
$n\in\dlat\setminus\{0\}$
either (i) of (ii) is true but the restriction $f|_I$ is discontinuous at some $\xi_1\in I$.
By Lemma \ref{lem:disco} there exists $n_0\in\dlat\setminus\{0\}$, such that $\xi_1+ n_0\in\cA$ and
$|g(\xi_1 +n_0)-g(\xi_1)|<2\rho^{-N}$. It follows in particular that $I+n_0\subset\cA$ 
hence (ii) above is true by our assumptions.

We now apply Lemma \ref{lem:athglob} with $a=\xi_1$ and $b=\xi_0$ (we may assume that $\alpha>2l-3$).
We conclude that there exist $m_1,m_2\in\dlat$, $m_1\neq m_2$, such that
\[
|g(\xi_0+m_1)-g(\xi_1)|\leq c(\rho^{2l-1}T +\rho^{-N}) \;\; , \;\;
|g(\xi_0+m_2)-g(\xi_1+n_0)|\leq c(\rho^{2l-1}T +\rho^{-N}).
\]
At least one of the $m_1,m_2$ is non-zero. If $m_1\neq 0$ then, using also estimate
(ii) of Proposition \ref{lem:illi}, we have
\begin{eqnarray*}
| g(\xi_0+m_1)-g(\xi_0)|&\leq& | g(\xi_0+m_1)-g(\xi_1)|+|g(\xi_1) -g(\xi_0)| \\
&\leq& c(\rho^{2l-1}T +\rho^{-N}) + c\rho^{2l-1}T \\
&\leq & c(\rho^{2l-1}T +\rho^{-N}).
\end{eqnarray*}
This contradicts (ii) (provided $c_2$ has been chosen to be large enough). Suppose now that $m_2\neq 0$. Then
\begin{eqnarray*}
| g(\xi_0+m_2)-g(\xi_0)|&\leq& | g(\xi_0+m_2)-g(\xi_1+n_0)|+|g(\xi_1+n_0) -g(\xi_1)| + \\
&&\quad + |g(\xi_1)-g(\xi_0)| \\
&\leq&  c(\rho^{2l-1}T +\rho^{-N}) + c\rho^{-N} + c\rho^{2l-1}T \\
&\leq & c(\rho^{2l-1}T +\rho^{-N}),
\end{eqnarray*}
which again contradicts (ii). This completes the proof of the lemma. $\hfill \Box$

\

Let us now write $\cD(\delta)$ as a disjoint union,
$\cD(\delta)=\cD_0(\delta)\cup\cD_1(\delta)\cup\cD_2(\delta)$, where: $\cD_0(\delta)$ contains
all $\xi\in\cD(\delta)$ for which there does not exist any $n\in\dlat\setminus\{0\}$ such that $\xi-n\in\cB(\delta)$;
$\cD_1(\delta)$ contains those $\xi\in\cD(\delta)$ for which there exists exactly one such $n$; and
$\cD_2(\delta)$ contains those $\xi\in\cD(\delta)$ for which there exist at least two (different) such points $n$.
\begin{lemma}
Suppose $d\geq 3$. Then there holds
\begin{eqnarray*}
\ia && \cB(\delta) \bigcap \Big(\bigcup_{n\in\dlat\setminus\{0\}} (\cD_0(\delta)-n) \Big)=\emptyset \; ; \\
\ib && \bigcup_{n\in\dlat\setminus\{0\}}(\cD_2(\delta)-n)\subset\bigcup_{n\in\dlat\setminus\{0\}}(\cB(\delta)-n) \; ; \\
\ic && \vol\Big( \bigcup_{n\in\dlat\setminus\{0\}}\Big( (\cB(\delta)-n) \cap \cB(\delta) \Big)\Big) \leq
c (\delta^2\rho^{-4l+2d+1} +\delta\rho^{(\alpha-2l+1)d-\alpha +d-1} ) \; ; \\
\id && \vol\Big( \bigcup_{n\in\dlat\setminus\{0\}}\Big( (\cD_1(\delta)-n) \cap \cB(\delta) \Big)\Big)\leq \delta\rho^{d-2l-\epsilon_0}\; .
\end{eqnarray*}
\label{lem:mar}
If $d=2$ then the same estimates are valid provided (iii) is replaced by
\[
{\rm (iii')} \qquad \vol\Big( \bigcup_{n\in\dlat\setminus\{0\}}\Big( (\cB(\delta)-n) \cap \cB(\delta) \Big)\Big) \leq
c (\delta^{3/2}\rho^{5-3l} +\delta\rho^{\alpha -4l+3} ) \; . \\
\]
\end{lemma}
{\em Proof.} Parts (i) and (ii) follow easily from the definition of the
sets $\cD_0(\delta)$ and $\cD_2(\delta)$. To prove (iii) we first note that from Lemma \ref{lem:volinters}
we have that for any $n\in\dlat\setminus\{0\}$ there holds
\begin{equation}
\vol ((\cB(\delta)-n)\cap \cB(\delta)) \leq c (\delta^2\rho^{-4l+d+1} + \delta\rho^{(\alpha-2l+1)d-\alpha} );
\label{win1}
\end{equation}
moreover, if in addition $n\in\dlat\setminus\{0\}$ satisfies $| |n|-2\rho|\geq c$, then we can do better, namely
\begin{equation}
\vol ((\cB(\delta)-n)\cap \cB(\delta)) \leq c\delta^2\rho^{-4l+d+1}.
\label{win11}
\end{equation}
But the number of $n\in\dlat\setminus\{0\}$ for which $| |n|-2\rho|\leq c$ 
is smaller than
$c\rho^{d-1}$; and the number  of $n\in\dlat\setminus\{0\}$ 
for which $(\cB(\delta)-n)\cap \cB(\delta)$ is non-empty is smaller than $c\rho^d$.
Hence (iii) follows. Finally, we have
\begin{eqnarray*}
&&\hspace{-2cm}\vol\Big( \bigcup_{n\in\dlat\setminus\{0\}}\Big( (\cD_1(\delta)-n) \cap \cB(\delta) \Big)\Big) \\
&\leq& \sum_{n\in\dlat\setminus\{0\}}\vol(\{\xi\in\cB(\delta) \; : \; \xi+n\in\cD_1(\delta)\})\\
&=&\sum_{n\in\dlat\setminus\{0\}}\vol(\{\eta\in\cD_1(\delta) \; : \; \eta -n\in\cB(\delta)\})\\
&=& \vol(\{\eta\in\cD_1(\delta_1) \; : \; \eta -n\in\cB(\delta)\, , \mbox{ for some }n\in\dlat\setminus\{0\}\})\\
&=&\vol(\cD_1(\delta))\\
&\leq&\rho^{d-2l-\epsilon_0} \, ,
\end{eqnarray*}
which is (iv). $\hfill \Box$

It is now the time to choose precise values of $N$ and $\delta$. We put
$\delta=c_3\rho^{2l-d-1}$, where $c_3$ is a (small) constant to be determined later, 
and $N=d+2$, so that for large $\rho$ we have
\begin{equation}
2\rho^{-N}\leq\delta \; .
\label{simple}
\end{equation}
For any unit vector $\eta\in\R^d$ let us define $I_{\eta}=\{ r\eta \; : \; r>0 \; , \;\; r\eta \in\cA(\delta)\}$.
We then have
\begin{lemma}
Let $\delta=c_3\rho^{2l-d-1}$ if $d\geq 3$ and $\delta=c_3\rho^{2l-6}$ if $d=2$.
If $c_3$ is small enough
then there exists at least one $\eta\in\R^d$, $|\eta|=1$, such that $I_{\eta}\subset\cB$ and
the restriction $f|_{I_{\eta}}$ is continuous.
\label{lem:ieta}
\end{lemma}
{\em Proof.} Let us assume the contrary. Lemma \ref{lem:c3} then implies that for any such interval $I_{\eta}$
and for any $\xi\in I_{\eta}$ there exists an $n\in\dlat\setminus\{0\}$ such that
\begin{equation}
\xi+n\in\cA  \quad , \quad |g(\xi+n)-g(\xi)|\leq c_2(T\rho^{2l-1}+\rho^{-N}).
\label{st}
\end{equation}
Since all such intervals $I_{\eta}$ cover $\cB(\delta)$, the existence
of an $n\in\dlat\setminus\{0\}$ satisfying (\ref{st})
is in fact true for any $\xi\in\cB(\delta)$. But (cf. (\ref{xt})) the length
of each such interval $I_{\eta}$ is $\asymp \delta\rho^{-2l+1}$, hence (\ref{st}) gives
$|g(\xi+n)-g(\xi)|\leq c(\delta+\rho^{-N})$. It follows that
\begin{eqnarray*}
|g(\xi+n)-\rho^{2l}|&\leq& |g(\xi+n)-g(\xi)| +|g(\xi)-\rho^{2l}| \\
&\leq& c(\delta+\rho^{-N}) + \delta \\
&\leq& (1+c)\delta \\
&=:&\delta_1\, .
\end{eqnarray*}
Hence we have proved that for any $\xi\in\cB(\delta)$ there exists $n\in\dlat\setminus\{0\}$ such that
$\xi+n\in\cA(\delta_1)$, that is $\cB(\delta)\subset\bigcup_{n\in\dlat\setminus\{0\}}(\cA(\delta_1)-n)$.
Recalling that $\cA(\delta_1)=\cB(\delta_1)\cup \cD(\delta_1)$ and using (i) and (ii) of Lemma \ref{lem:mar} we
obtain
\begin{equation}
\cB(\delta)\subset\bigcup_{n\in\dlat\setminus\{0\}}(\cB(\delta_1)-n)
\bigcup\bigcup_{n\in\dlat\setminus\{0\}}(\cD(\delta_1)-n).
\label{gl}
\end{equation}
Combining this with (i) and (ii) of Lemma \ref{lem:mar} gives
\begin{equation}
\cB(\delta)=\bigcup_{n\in\dlat\setminus\{0\}}\Big((\cB(\delta_1)-n)\cap\cB(\delta)\Big)
\bigcup\bigcup_{n\in\dlat\setminus\{0\}}\Big((\cD_1(\delta_1)-n)\cap\cB(\delta)\Big).
\label{gl1}
\end{equation}
We now consider the respective volumes in (\ref{gl1}). Assume that $d\geq 3$.
Using part (ii) of Lemma \ref{lem:volumes} and parts (iii) and (iv) of Lemma \ref{lem:mar} we conclude that
\begin{eqnarray*}
\delta\rho^{d-2l}&\leq&c(\delta_1^2\rho^{-4l+2d+1} +\delta_1\rho^{(\alpha-2l+1)d-\alpha +d-1})+\delta_1 \rho^{d-2l-\epsilon_0}\\
&\leq& c(\delta^2\rho^{-4l+2d+1}+ \delta\rho^{(\alpha-2l+1)d-\alpha +d-1}+\delta\rho^{d-2l-\epsilon_0}) \\
&=:& I_1+I_2+I_3 \; .
\end{eqnarray*}
We have $I_3=o(\delta\rho^{d-2l})$. Since $d>1$ we also have $I_2=o(\delta\rho^{d-2l})$. Hence we conclude that
\[
\delta\rho^{d-2l}\leq c\delta^2\rho^{-4l+2d+1}.
\]
Recalling that $\delta=c_3\rho^{2l-d-1}$, we reach a contradiction if $c_3$ is small
enough. The same argument works when $d=2$. $\hfill \Box$

\begin{theorem}
Let $d\geq 3$. Suppose that $\rho$ is large enough. Then $\lambda=\rho^{2l}$ belongs in
$\sigma(H)$. Moreover, there exists $c_3>0$ such that the interval
$[\rho^{2l}-c_3\rho^{2l-d-1},\rho^{2l}+c_3\rho^{2l-d-1}]$ lies inside a single spectral band.
If $d=2$ then the same is true, but the respective interval is $[\rho^{2l}-c_3\rho^{2l-6},\rho^{2l}+c_3\rho^{2l-6}]$.
\label{thm:rfin}
\end{theorem}
{\em Proof.} Let $d\geq 3$. We may assume that $\alpha\geq 2l-d-1$.
Let $I_{\eta}$ be an interval with the properties specified in Lemma \ref{lem:ieta}. Then the value of $f$ at the one end of $I_{\eta}$
is $\rho^{2l}+c_3\rho^{2l-d-1}$, and the value at the other end is $\rho^{2l}-c_3\rho^{2l-d-1}$. Since $f|_{I_{\eta}}$ is continuous, it takes the value $\rho^{2l}$. When $d=2$ we argue similarly. $\hfill \Box$

\

\parindent=0pt
\parskip=12pt

{\bf Acknowledgment.} The first author acknowledges the warm hospitality of the Department of Mathematics
of University College London where part of this work was carried out.
The work of the second author 
was partially supported by the Leverhulme trust and by the EPSRC grant EP/F029721/1.

%%%%%%%%%%%%%%%%%%%%%%%%%%%%%%%%%%%%%%%%%%%%%%%%%%%%%%%%%

G. Barbatis: Department of Mathematics, University of Athens, Panepistimioupolis, 15784 Athens, Greece

L. Parnovski: Department of Mathematics, University College London, Gower street, London WC1E 6BT, UK

\end{document}